\documentclass{amsart}

\usepackage{epsfig,amssymb,amsmath,float,subfig}

\newcommand{\real}{{\mathbb{R}}}

\newcommand{\sta}{{\rm sta}}
\newcommand{\ext}{{\rm ext}}

\newcommand{\cal}{\mathcal}

\newcommand{ \calI }{{\cal{I}}}

\newcommand{ \PP}{{P}}

\newcommand{\uu}{u}

\newcommand{\half}{\frac{1}{2}}

\newcommand{\bsig}{\mbox{\boldmath$\sigma$}}
\newcommand{\bvsig}{\mbox{\boldmath$\varsigma$}}

\newcommand{\bveps}{\mbox{\boldmath$\varepsilon$}}

\newcommand{\barbvsig}{\bar{\mbox{\boldmath$\varsigma$}}}

\newcommand{\beps}{\mbox{\boldmath$\epsilon$}}

\newcommand{\btau}{{\mbox{\boldmath$\tau$}}}

\newcommand{\bgamma}{{\mbox{\boldmath$\gamma$}}}

\newcommand{\la}{\langle}
\newcommand{\ra}{\rangle}

\newcommand{\ba}{{{\bf a}}}

\newcommand{\Oo}{\Omega}

\newcommand{\bE}{{\bf E}}
\newcommand{\bQ}{{\bf Q}}

\newcommand{\bM}{{\bf M}}
\newcommand{\bK}{{\bf K}}
\newcommand{\eba}{\begin{array}}
\newcommand{\eea}{\end{array}}
\newcommand{\ebe}{\begin{eqnarray}}
\newcommand{\eee}{\end{eqnarray}}
\newcommand{\eb}{\begin{equation}}
\newcommand{\ee}{\end{equation}}
\newcommand{\bT}{{\bf T}}
\newcommand{\bJ}{{\bf J}}

\newcommand{\calP}{{\cal{P}}}

\newcommand{\bC}{{\bf C}}
\newcommand{\bG}{{\bf G}}

\newcommand{\bff}{{\bf f}}

\newcommand{\bR}{{\bf R}}

\newcommand{\bx}{{\bf x}}

\newcommand{\bN}{{\bf N}}
\newcommand{\bF}{{\bf F}}
\newcommand{\bA}{{\bf A}}
\newcommand{\bB}{{\bf B}}
\newcommand{\bI}{{\bf I}}

\newcommand{\calS}{{\cal S}}

\newcommand{\calSa}{\calS_a}

\newcommand{\calU}{{\cal U}}

\newcommand{\calV}{{\cal V}}

\newcommand{\bu}{{\bf u}}

\newcommand{\barbx}{\bar{\bf x}}

\newcommand{\bgra}{{\bf \nabla}}

\newcommand{\bP}{{\bf P}}

\newcommand{\dO}{\,\mbox{d}\Oo}

\newcommand{\bU}{{\bf U}}

\newcommand{\alp}{{\alpha}}

\newcommand{ \eps}{{\epsilon}}

\newcommand{ \sig}{{\sigma}}
\newcommand{ \Lam}{{\Lambda}}

\newcommand{ \lam}{{\lambda}}

\newcommand{\bD}{{\bf D}}

\newcommand{\tr}{{\mbox{tr}}}

\newcommand{\Rank}{{\mbox{Rank }}}
\newcommand{\Diag}{{\mbox{Diag }}}

\newtheorem{lemma}{Lemma}
\newtheorem{remark}{Remark}

\newtheorem{thm}{Theorem}
\renewcommand\bI{{{I}}}
\renewcommand\barbx{{\bar{x}}}

\renewcommand\beta{{{b}}}

   \renewcommand\Im{I_m}
      \newcommand\Ip{I_p}

\renewcommand\eb{\begin{equation}}
\renewcommand\ee{\end{equation}}
\renewcommand\bar{\overline}
\newcommand\Col{{\cal C}_{ol}}
    \newcommand\WW{W}
    \newcommand\VV{V}

\begin{document}

\title[Canonical Duality-Triality Theory in complex systems]{Canonical Duality-Triality Theory 
for solving general global optimization problems\\ in complex systems}


\author[D. Morales-Silva]{Daniel Morales-Silva}
\address{School of Science, Information Technology and Engineering, Federation University Australia, Victoria 3353, Australia}
\email{d.moralessilva@federation.edu.au}

\author[D. Y. Gao]{David Y. Gao}
\address{School of Science, Information Technology and Engineering, Federation University Australia, Victoria 3353, Australia\\
Research School of Engineering, Australian National University, Canberra, Australia}
\email{d.gao@federation.edu.au, david.gao@anu.edu.au}



\begin{abstract}
General nonconvex optimization problems are studied by using  the canonical duality-triality theory. The
triality theory is proved for sums of exponentials and quartic polynomials, which solved
an open problem left in 2003. This theory   can be used to find the global minimum and local extrema,
which  bridges a  gap between global  optimization and nonconvex mechanics.
Detailed applications are illustrated
  by several  examples.
\end{abstract}


\maketitle
\section{Introduction and Motivation}
This  paper intends to  solve the following nonconvex optimization problem ($(\calP)$ in short):
\begin{equation}\label{ProbP}
(\calP):  \;\;\; \ext\left\{\Pi(\bx)=W(\bx)+\frac{1}{2}\bx^t\bA\bx-{\bf
f}^t\bx \; | \; \bx \in \real^n \right\},
\end{equation}
where $\ext \{* \} $ denotes finding extremum points of a function given in $\{ * \}$,
  ${\bf f}\in\real^n$ is a given (input) vector,  $\bA\in\real^{n\times n}$ is a given symmetric matrix, and $W:\real^n\rightarrow \real$ is a combination of fourth order polynomials (double-well functions) and quadratic-exponential functions, namely:
\[
W(\bx):=\sum_{i\in\Im} \exp\left(\frac{1}{2}\bx^t\bB_i\bx-\alp_i\right)+
\sum_{j\in\Ip}\frac{\beta_j}{2}\left(\frac{1}{2}\bx^t\bC_j\bx-\theta_j\right)^2,
\]
where $\Im = \{ 1 , \dots , m\}, \; \Ip=\{ 1, \dots, p\}$ are two integer sets with
$m, p$ which are  fixed integers; all the coefficients $\beta_j$ with $j\in I_p$ are positive constants, and $\alp_i , \theta_j \in \real \;\; \forall i \in \Im,  \; j \in \Ip$ are given parameters;  the matrices $\displaystyle\{\bB_i\}_{i\in \Im} $ and $ \{\bC_j\}_{j\in \Ip}$ are assumed to be symmetric, positive semi-definite such that the cone generated by them contains a positive definite matrix.\\

  The nonconvex optimization
problem $(\calP)$ arises  naturally in complex systems with  a wide range of  applications, including
chaotical dynamical systems \cite{Gao-amma03,DGAO2OG,DGAO3RU}, computational biology \cite{zgy},
chemical database analysis \cite{xie},
 large deformation  computational mechanics \cite{gao-jem,santos-gao},  population growing \cite{ruan-gao-ima},
 location/allocation,
network communication \cite{gao-ruan-pardalos}, and
   phase transitions
of solids \cite{DGAO2OG,gao-ogden08b,gao-yu}, etc.

For example, the popular sensor   network location problem
  is to  solve the following   system of
nonlinear  equations (see \cite{asp-god,more-wu}):
\begin{eqnarray}
\|\bu_i- \bu_j\|_2^2 = d_{ij}^2,\  \; \forall (i , j) \in \calI_p, \;\;
  \bu_k =  \ba_k ,\; \;\; \forall k \in \calI_b  \label{s00}
\end{eqnarray}
 where the vectors  $\bu_i = \{ \uu^\alp_i \} \in \real^d$ ($i=1, \dots, p$)
 represent the locations of the unknown sensors,
$\calI_p =\{ (i,j): \; i < j, \; d_{ij} \mbox{ is specified} \} $
   and $\calI_b=\{ k :  \; \bu_k = \ba_k \mbox{ is specified } \}$
are two given index sets,
 $d_{ij} $ are given distances for $ (i,  j ) \in \calI_p$, the given vectors
 $\ba_1, \ba_2,\cdots, \ba_q \in \real^d$
 are the so-called anchors. The notation $\| \bu_i  - \bu_j\|_2$ denotes the Euclidian distance between $\bu_i$ and $\bu_j$,
 i.e.,
 \[
\| \bu_i  - \bu_j\|_2 = \sqrt{ \sum_{\alp = 1}^d ( \uu_i^\alp - \uu_j^\alp)^2 }.
\]
By using the least squares method,
the quadratic equations (\ref{s00}) of the sensor
localization problem can be  reformulated as
an
optimization problem: 
\eb
\displaystyle
\min \left\{  \PP(\bu)= \displaystyle
\sum_{(i,j)\in \calI_p}
\frac{1}{2}\left( \| \bu_i- \bu_j\|_2^2 - d_{ij}^2 \right)^2  : \;\;
\bu_i \in \calU_a  \right\} , \label{eq-leasq}
\ee
where
$\calU_a = \{ \bu \in \real^{d \times p} | \;\; \bu_k = \ba_k \;\; \forall
k \in \calI_b\}$ is a feasible space.
Let $\bx = \{ \{ \uu_1^1, \dots, \uu^d_1\}, \dots, \{ \uu_p^1, \dots, \uu^d_p\} \} \in \real^n$ $(n=d \times p$)
denote an extended vector. By using
Lagrange multiplier method to relax the
boundary conditions in $\calU_a$, the least squares method for
the  sensor localization problem
(\ref{eq-leasq}) can be written in the problem (\ref{ProbP})
for certain properly defined matrices $\{ {\bf C}_j\}$,
which is the so-called deformation matrix in structural mechanics.
The sensor network localization type problems also appear
 in  computational biology, Euclidean ball packing,  molecular confirmation,
 and recently, wireless network communication, etc \cite{ruan-gao-pe,zgy}. Due to the nonconvexity,
   the sensor network localization problem is  considered  to be
   {  NP-hard} even for the simplest case $d=1$ \cite{more-wu,saxe}.
Recent result of Aspnes {\em et al} \cite{asp-god}
shows that the problem of computing a realization of the sensors on the plane is
NP-complete in general.

Mathematics and mechanics have been  two  complementary partners since the Newton times.
Many fundamental ideas, concepts, and mathematical methods extensively used in
calculus of variations and optimization are originated from mechanics.
For examples, the Lagrange multiplier method was first proposed by Lagrange from
the classical analytic mechanics; while the concepts of super-potential and
sub-differential in modern convex analysis were introduced by
Moreau from frictional mechanics \cite{moreau,more-pana-stra-88}.
From the point view of computational large deformation mechanics,
both the  fourth-order polynomial minimization problem $(\calP)$ and
the sensor localization problem (\ref{eq-leasq}) are
actually  two special cases of
   discretized finite deformation problems
 \cite{gao-jem}.
  It is known that in continuum mechanics and differential geometry, the
 deformation   $\bu(\bx):\Oo\rightarrow \real^r  $
 is  a vector field over an open domain $\Oo \subset\real^r$,
 and the minimal potential variational problem is defined by
  \eb
 \min \;\; \left\{ \PP(\bu) = \int_\Oo [  W(\nabla \bu) - \bu^T \bff ] \dO
   \; | \;\; \bu \in \calU_a \right\}, \label{eq-minpot}
  \ee
  where   $ W(\bF)$ is the so-called {\em stored strain energy}, which is usually a
 nonconvex function of the deformation gradient $\bF = \nabla \bu$,
 the feasible set  $\calU_a$ in this nonconvex variational problem is called
 the {\em kinematically admissible space}, where  certain boundary conditions
 are prescribed.
 According to the hyper-elasticity law (see Chapter 6.1.2 \cite{DGAO1} or \cite{marsden-hugh}),
 the stored strain energy should be an {\em objective function } of the deformation gradient $\bF$,
    i.e., there exists an objective strain measure $\bE (\bF)$ and a convex
  function  $\VV(\bE)$ such that
 \eb
 \WW(\nabla \bu) = \VV(\bE(\nabla \bu)) .
 \ee
 One of the most simple objective strain measures is the
 well-known {\em Green-St. Venant strain tensor} $\bE = \half [\bF^T \bF - \bI]$.
 Clearly, this strain measure satisfies the objectivity condition, i.e.
   $\bE ({\bf Q} \bF) = \bE(\bF)$ for any given  orthonormal  (rotation) matrix ${\bf Q}$.
 For the most simple {\em St. Venant-Kirchhoff material},  $\VV(\bE)$ is a
 quadratic function of $\bE$, i.e.
 \eb
 \VV(\bE) = \half \lambda \left( \tr \bE  \right)^2 + \mu \tr \left(\bE \right)^2,
 \ee
 where, $\lambda , \mu > 0$ are the classical Lam\'{e} constants, $\tr \bE$ represents the trace of $\bE$.
  Therefore, the stored energy
 $\WW(\bF)$ is a fourth-order polynomial tensor function of $\bF = \nabla \bu $.
 While for bio-materials, the stored energy could be the combination of the
 polynomial and exponential functions of the Cauchy-Green
 strain tensor.
 By using   finite difference method (FDM), the deformation gradient
 $\nabla \bu$ can be  directly
 approximated by the difference $\bD \bu = \bu (\bx_i) - \bu(\bx_j) = \bu_i - \bu_j$.
 While in  finite element method (FEM), the domain
 $\Oo = \bigcup_e^m  \Oo^e $ is discretized by
 a finite number of elements $\Oo^e \subset \Oo$ and in each element, the deformation field
 $\bu(x) = \sum_i \bN_i(\bx) \bu_i $ is numerically  represented by the nodal vectors $\bu_i$
 via piecewise   interpolation (polynomial) function $\bN_i(\bx)$ (cf. \cite{gao-jem}).
 Therefore, by either FDM or FEM, the
 minimal potential variational problem (\ref{eq-minpot}) can be  eventually reduced to a
 very complicated large-scale fourth-order polynomial/exponential minimization problem with
 the problems $(\calP)$  as its  the most simple case.
 In the contact mechanics and elasto-plastic design of large deformed structures,
 the nonconvex problems are usually subjected to
 inequality constraints. In these cases, the global optimal solution could be local minima
 (see \cite{cai-gao-qin}) and
  to solve such problems is fundamentally difficult by using traditional direct methods.

Canonical duality theory was developed originally from  Gao and Strang's work in 1989 \cite{gao-gs1}
for solving general
 variational problem  (\ref{eq-minpot}) in finite deformation theory, where the stored energy
 $W(\bF)$ is  nonconvex and even nonsmooth. By introducing a so-called complementary gap function,
 they recovered the complementary energy principle in large deformation (geometrically nonlinear)  systems.
  They proved that the nonnegative gap function can be used to identify the global minimizer of
  the nonconvex potential variational problems.
  Seven years later, it was discovered that the negative gap function can be used to identify the
largest local minimum and maximum. Therefore, a so-called
triality theory was first proposed in nonconvex mechanics \cite{gao-amr}, and then generalized to
global optimization
\cite{gao-jogo00}. This triality theory is composed of a canonical min-max duality and two pairs
of double-min, double-max dualities, which reveals an intrinsic duality pattern in
complex systems and has been used successfully for solving a
wide class of challenging problems in complex systems \cite{gao-ima, gao-mecc,gao-cace,gao-sherali09}. However, it
was realized in 2003 \cite{Gao-amma03, gao-opt03} that the double-min duality holds under
"certain additional conditions". Recently, this problem is partly solved for a class
of fourth order polynomial optimization problems \cite{WU}.
   Based on these results, this paper  intends to solve the
more challenging problem $(\calP)$.
 We will show that by the canonical dual transformation,
  all  critical solutions of  $(\calP)$ can be analytically
presented in terms of the canonical dual solutions.
The extremality of these solutions  can be identified by the triality theory.
Several solved examples are listed in the last section.\\

\section{Canonical Dual Problem and Analytical Solutions}

Following the standard procedure of the canonical dual transformation
(cf. e.g., \cite{gao-opt03}), first we need to choose a geometric operator
$\Lam=(\Lam_1(\bx),\Lam_2(\bx)):\real^n\rightarrow \real^{m+p}$, where
\begin{eqnarray*}
 \Lam_1(\bx) &=& \left\{ \frac{1}{2}\bx^t\bB_i\bx-\alp_i \right\}  : \; \real^n \rightarrow \real^m, \\
 \Lam_2(\bx) &=& \left\{ \frac{1}{2}\bx^t\bC_j\bx-\theta_j \right\} : \real^n \rightarrow \real^p.
\end{eqnarray*}
Therefore, the nonconvex function $W(\bx)$ can be written in the following canonical form
\eb
W(\bx) = V(\Lam(\bx)) = V_1 (\Lam_1(\bx)) + V_2(\Lam_2(\bx))
\ee
with
\eb
 V_1(\beps)=\sum_{i\in\Im} \exp(\eps_i), \;\; V_2(\bgamma) = \sum_{j\in \Ip} \half \beta_j \gamma_j^2 .
\ee
Clearly, the canonical function $V(\bveps)$ is convex on
\eb
\calV_a = \{ \bveps = ( \beps, \bgamma) \in \real^{m+p} |
\;\;\eps_i \in [ -\alp_i, +\infty), \;\; \gamma_j \in [-\theta_j, +\infty), \;\; \forall i\in \Im, \; j\in \Ip\}
\ee
such that the canonical dual variable $\bvsig = (\btau, \bsig)$ of $\bveps = ( \beps, \bgamma)$
can be uniquely defined by
\eb
\bvsig =  \nabla V(\bveps) \Rightarrow  \;\;  \btau = \nabla V_1(\beps) = \{ \exp(\eps_i) \} , \;\;
\bsig = \nabla V_2(\bgamma) = \{ \beta_j \gamma_j \},
\ee
and on the canonical dual space
\eb
\calV^*_a = \{ \bvsig = ( \btau, \bsig) \in \real^{m+p} | \; \tau_i \in [ \exp(-\alp_i), \infty),
\;\; \sig_j \in [-\beta_j \theta_j , \infty), \; \forall i\in\Im, \; j\in\Ip \},
\ee
the Legendre conjugate of $V(\bveps)$ can be defined by
\eb
V^c(\bvsig) =  \sta\{ \bveps^t \bvsig - V(\bveps ) | \; \bveps \in \calV_a \}
= V^c_1(\btau) + V^c_2(\bsig)
\ee
where $\sta\{ *\}$ denotes finding stationary points of the function given in $\{ *\}$ and
\[
V^c_1(\btau) =  \sum_{i\in I_m}\left(\tau_i\ln\tau_i-\tau_i\right) , \;\; V^c_2(\bsig) =
 \sum_{j\in I_p}\frac{1}{2\beta_j}\sig^2_{j}.
 \]
  By using the canonical dual transformation
  $W(\bx) = V(\Lam(\bx)) = \Lam(\bx)^T \bvsig - V^c(\bvsig)$,
  the Gao-Strang {\em total complementary function}
  $\Xi:\real^n\times \calV_a^*\rightarrow \real$ associated with  the problem $(\calP)$ can be given by
  \begin{eqnarray}
  \Xi(\bx,\bvsig) & = & \la
\Lam (\bx),\bvsig\ra - V^c(\bvsig)+\frac{1}{2}\bx^t\bA\bx-{\bf f}^t\bx \nonumber \\
&= & \frac{1}{2}\bx^t\bG(\bvsig)\bx- {\bf \alp}^t\btau-{\bf\theta}^t\bsig -V_1^c(\btau)-V_2^c(\bsig) -{\bf f}^t\bx,
 \label{xireduced}
\end{eqnarray}
where
\eb
\displaystyle \bG(\bvsig)=\bA+\sum_{i\in I_m}\tau_i\bB_i+\sum_{j\in I_p}\sig_j\bC_j.
 \ee
 Via this $\Xi(\bx,\bvsig)$, the canonical dual function
 $\Pi^d:\calV_a^*\rightarrow \real$ can be defined by
\[
\Pi^d(\bvsig):= \sta \left\{  \Xi(\bx,\bvsig) | \; \bx \in \real^n \right\}
 =\left\{\Xi\left(\bx(\bvsig),\bvsig\right) : \bgra_{\bx}\Xi\left(\bx(\bvsig),\bvsig\right)=0\right\}.
 \]

Notice that $\bgra_{\bx}\Xi(\bx,\bvsig)=\bG(\bvsig)\bx-{\bf f}=0 $
if and only if
\eb
\bG(\bvsig)\bx={\bf f}.\label{xeqn}
\ee
Let $\Col(\bG(\bvsig))$ be the space generated by the columns of the matrix $\bG(\bvsig)$.
Then, on the dual feasible space
\[
\calS_a=\left\{\bvsig\in\calV^*_a:{\bf f}\in \Col(\bG(\bvsig))\right\},
\]
the primal solution $\bx=(\bG(\bvsig))^{-1}{\bf f}$ is well defined (if $\bG(\bvsig)$ is singular, $(\bG(\bvsig))^{-1}$ denotes its pseudo-inverse, see \cite{DESOER}, \cite{PETERS} and references therein) and we have $\Pi^d:\calS_a\rightarrow \real$
\begin{eqnarray}\label{PidDefin}
\Pi^d(\bvsig) &=& -\frac{1}{2}{\bf f}^t(\bG(\bvsig))^{-1}{\bf f}-V_1^c(\btau)-V_2^c(\bsig)- {\bf \alp}^t\btau-{\bf \theta}^t\bsig.
 \end{eqnarray}
Therefore, the canonical dual problem is proposed in the following  form:
\begin{equation}\label{ProbPd}
(\calP^d): \;\; \ext \{\Pi^d(\bvsig):\bvsig\in\calS_a\}.\end{equation}

By the canonical duality theory, it is not difficult to show that
\begin{equation}\label{pieqxi}
\Pi(\bx)=\sta \{ \Xi(\bx,\bvsig) : \bvsig \in \calS_a \}
=\Xi(\bx,\bvsig(\bx)),\end{equation}
where
$\bvsig(\bx)=(\btau(\bx),\bsig(\bx))$ and
$$(\btau(\bx))_i= \exp((\Lam_1(\bx))_i),\ i\in I_m,$$
 $$(\bsig(\bx))_j= \beta_{j} (\Lam_2(\bx))_{j},\ j\in I_p.$$

According to the general theory presented in \cite{gao-opt03}, we have the following result.
\begin{thm}[Analytical Solutions]\label{thm-1}
Suppose that for a given ${\bf f} \in \real^n$ the canonical dual space $\calS_a$ is not empty.
If  $\overline{\bvsig}\in \calS_a$ is a stationary point of $\Pi^d$, then
\eb
\overline{\bx}=(\bG(\overline{\bvsig}))^{-1}{\bf f}
\ee
is a stationary point of $\Pi$ and
\eb
\Pi(\overline{\bx})=\Pi^d(\overline{\bvsig}).
\ee
\end{thm}
\textbf{Proof:}
Let us calculate $\bgra \Pi^d(\bvsig)$ and $\bgra^2 \Pi^d(\bvsig)$. We know that $$\bgra \Pi^d(\bvsig)=\left[\begin{array}{c}
                                                                 \bgra_{\btau}\Pi^d(\bvsig) \\
                                                                 \bgra_{\bsig}\Pi^d(\bvsig)
                                                               \end{array}\right] \in \real^{m+p}, $$ then \begin{equation}\label{gradPidp1} (\bgra_{\btau}\Pi^d(\bvsig))_i=\frac{1}{2}{\bf f}^t(\bG(\bvsig))^{-1}\bB_i(\bG(\bvsig))^{-1}{\bf f}-\ln\tau_i-\alp_i,\ i\in I_m; \end{equation}
\begin{equation}\label{gradPidp2} (\bgra_{\bsig}\Pi^d(\bvsig))_j=\frac{1}{2}{\bf f}^t(\bG(\bvsig))^{-1}\bC_j(\bG(\bvsig))^{-1}{\bf f}-\frac{\sig_j}{\beta_{j}}-\theta_j,\ j\in I_p.
\end{equation}

On the other hand,
$$ \bgra^2 \Pi^d(\bvsig)=\left[\begin{array}{cc}
\bgra^2_{\btau\btau}\Pi^d(\bvsig) & \bgra^2_{\btau\bsig}\Pi^d(\bvsig) \\
\bgra^2_{\bsig\btau}\Pi^d(\bvsig) & \bgra^2_{\bsig\bsig}\Pi^d(\bvsig)
                                                        \end{array}
\right] \in \real^{m+p}\times\real^{m+p},$$ where $\bgra^2_{\btau\bsig}\Pi^d(\bvsig):=(\bgra_{\btau}(\bgra_{\bsig}\Pi^d(\bvsig))^t)$. Let $\delta_{ij}$ be the Kronecker's delta. Then

\begin{eqnarray*}
(\bgra^2_{\btau\btau}\Pi^d(\bvsig))_{ij} & = & -{\bf f}^t(\bG(\bvsig))^{-1}\bB_i(\bG(\bvsig))^{-1}\bB_j(\bG(\bvsig))^{-1}{\bf f}-\frac{\delta_{ij}}{\tau_j},\\
& &  i,j\in I_m.  \\
(\bgra^2_{\btau\bsig}\Pi^d(\bvsig))_{ij} & = & -{\bf f}^t(\bG(\bvsig))^{-1}\bB_i(\bG(\bvsig))^{-1}\bC_j(\bG(\bvsig))^{-1}{\bf f} \\
& & i\in I_m; j\in I_p.\\
(\bgra^2_{\bsig\btau}\Pi^d(\bvsig))_{ij} & = & -{\bf f}^t(\bG(\bvsig))^{-1}\bC_i(\bG(\bvsig))^{-1}\bB_j(\bG(\bvsig))^{-1}{\bf f} \\
& & i\in I_m; j\in I_p. \\
(\bgra^2_{\bsig\bsig}\Pi^d(\bvsig))_{ij} & = & -{\bf f}^t(\bG(\bvsig))^{-1}\bC_i(\bG(\bvsig))^{-1}\bC_j(\bG(\bvsig))^{-1}\bff-\frac{\delta_{ij}}{\beta_{j}}\\
& & i,j\in I_p.
\end{eqnarray*}

By making $\bx=(\bG(\bvsig))^{-1}{\bf f}$ and $\bF(\bx)\in \real^{n\times(m+p)}$ be \newline $\bF(\bx)=[\bB_1\bx,\ldots,\bB_{m}\bx,\bC_1\bx,\ldots,\bC_p\bx]$, we have
\begin{equation}\bgra^2\Pi^d(\bvsig)=-\bF(\bx)^t(\bG(\bvsig))^{-1}\bF(\bx)-
\Diag\left(\frac{1}{\tau_1},\ldots,\frac{1}{\tau_{m}},\frac{1}{\beta_1},\ldots,\frac{1}{\beta_{p}}\right).\end{equation}

Let $\displaystyle\bD=\Diag\left(\tau_1,\ldots,\tau_{m},\beta_1,\ldots,\beta_{p}\right)$, then $\bgra^2\Pi^d(\bvsig)$ can be written as

\begin{equation}\label{HessianPid}\bgra^2\Pi^d(\bvsig)=-\bF(\bx)^t(\bG(\bvsig))^{-1}\bF(\bx)-\bD^{-1}.\end{equation}

Calculating $\bgra \Pi(\bx)$ and $\bgra^2\Pi(\bx)$, we have respectively
\begin{equation}\label{gradPi}\bgra\Pi(\bx)=\sum_{i\in I_m}\exp\left(\frac{1}{2}\bx^t\bB_i\bx-\alp_i\right)\bB_i\bx+\sum_{j\in I_p}\beta_j\left(\frac{1}{2}\bx^t\bC_{j}\bx-\theta_{j}\right)\bC_{j}\bx+\bA\bx-{\bf f}.\end{equation}
\begin{eqnarray}\label{HessianPi}
\bgra^2\Pi(\bx)& = &  \bA + \sum_{i\in I_m}\exp\left(\frac{1}{2}\bx^t\bB_i\bx-\alp_i\right)(\bB_i\bx(\bB_i\bx)^t+\bB_i)\nonumber \\
&  & +\sum_{j\in I_p}\beta_j\left(\bC_{j}\bx(\bC_{j}\bx)^t
+\left(\frac{1}{2}\bx^t\bC_{j}\bx-\theta_{j}\right)\bC_{j}\right)
.
\end{eqnarray}

Since $\overline{\bvsig}=(\overline{\btau},\overline{\bsig})$ is a stationary point of $\Pi^d$ then by Equations \eqref{gradPidp1} and \eqref{gradPidp2} we have that \begin{equation}\label{sigmaandlambdap1} (\Lam_1(\overline{\bx}))_i=\ln\overline{\tau}_i,\ i\in I_m ;\end{equation}
\begin{equation}\label{sigmaandlambdap2} (\Lam_2(\overline{\bx}))_j=\frac{\overline{\sig}_j}{\beta_{j}},\ j\in I_p. \end{equation}
Using Equations \eqref{sigmaandlambdap1} and \eqref{sigmaandlambdap2} in Equation \eqref{gradPi}, we obtain $$\bgra\Pi(\overline{\bx})=\bG(\overline{\bvsig})\overline{\bx}-{\bf f}=\bG(\overline{\bvsig})(\bG(\overline{\bvsig}))^{-1}{\bf f}-{\bf f}=0.$$
Notice that Equations \eqref{sigmaandlambdap1} and \eqref{sigmaandlambdap2} together with Equations \eqref{PidDefin} and \eqref{pieqxi} imply that \begin{equation}\label{pixeqspidsig}\Pi(\overline{\bx})= \Xi(\overline{\bx},\overline{\bvsig})=\Xi((\bG(\overline{\bvsig}))^{-1}{\bf f},\overline{\bvsig})=\Pi^d(\overline{\bvsig}).\end{equation} And
this finishes the proof.\hfill $\blacksquare$

\begin{remark}
This theorem shows that the problem $(\calP^d)$ is canonical dual to the nonconvex primal problem $(\calP)$
in the sense that $\Pi(\barbx) = \Pi^d(\barbvsig)$ at each critical point of $\Xi(\bx, \bvsig)$.
By the criticality condition \eqref{xeqn} we know that if $\bG(\bvsig)$ is singular at $\barbvsig$, the
canonical equilibrium equation  \eqref{xeqn} may have infinite number of solutions:
$\barbx = \bG(\barbvsig)^{\dag}  \bff + \bN \bx^o$, where $\bG^{\dag}$ represents the Moore-Penrose generalized inverse, $\bN$ is a basis matrix of the null space of $\bG(\barbvsig)$,  and $\bx^o$ is a free vector.
In this case, Theorem \ref{thm-1} still holds, but the canonical dual function $\Pi^d$ will have additional parametrical vector $\bx^o$. In order to avoid this case, a quadratic perturbation method is introduced
in \cite{ruan-gao-pe}, i.e. in the case that $\bG(\barbvsig)$ is singular, replace it by the following perturbed
 form
 \eb
 \bG_\alp(\barbvsig) = \bG(\barbvsig) + \alp \bD
 \ee
 where $\alp> 0 $ is a perturbation parameter and $\bD$ is a given positive-definite matrix. Very often, $\bD = {\bf I}$. Detailed study on this quadratic perturbation method is given in \cite{ruan-gao-pe}.
\end{remark}

In the next section, we will show that the extremality of some of these solutions can be identified by a refined triality theory.

\section{Triality Theory}

Before presenting the refined triality theory, we  need  the following sets
$$\calS_a^+:=\{\bvsig\in\calS_a:G(\bvsig)\succeq 0\},\quad \calS_a^-:=\{\bvsig\in\calS_a:G(\bvsig)\prec 0\}.$$

\begin{lemma}\label{MatrixL}
Suppose that $m+p<n$, $\overline{\bvsig}\in\calS_a^-$ is a stationary point and a local minimizer of $\Pi^d$ and $\bar{\bx}=(\bG(\overline{\bvsig}))^{-1}{\bf f}$. Then, there exists a matrix ${\bf L}\in\real^{n\times(m+p)}$ with $\Rank({\bf L})=m+p$ such that \begin{equation}\label{localminm1m2<n}{\bf L}^t\bgra^2\Pi(\overline{\bx}){\bf L}\succeq 0.\end{equation}
\end{lemma}
\textbf{Proof:} Since $\overline{\bvsig}\in\calS_a^-$ is a local minimizer of $\Pi^d$, we have that $\bgra^2\Pi^d(\overline{\bvsig})\succeq 0$.
It follows from Equation \eqref{HessianPid} that $$-\bF(\overline{\bx})^t(\bG(\overline{\bvsig}))^{-1}\bF(\overline{\bx})\succeq \bD^{-1}\succ 0.$$
Thus, $\Rank(\bF(\overline{\bx}))=m+p$. Since $\overline{\bvsig}\in\calS_a^-$ and $\bF(\overline{\bx})\bD\bF(\overline{\bx})^t\succeq 0$
there exists a nonsingular matrix $\bT\in\real^{n\times n}$ such that
\begin{equation}
\bT^t\bG(\overline{\bvsig})\bT=
\Diag(-\lam_1,\ldots,-\lam_n)\end{equation}
and
\begin{equation}\bT^t\bF(\overline{\bx})\bD\bF(\overline{\bx})^t\bT=\Diag(a_1,\ldots,a_{m_1+m_2},0,\ldots,0),\end{equation}
where $\lam_i>0$ for every $i=1,\ldots,n$ and $a_j>0$ for every $j=1,\ldots,m+p$
(see \cite{FENGLIN}, \cite{HORN} and references therein). According to Lemma \ref{SVDLEMMA} in the Appendix, we know that there exists orthogonal matrices $\bU\in\real^{n\times n}$ and
 $\bE\in\real^{(m+p)\times(m+p)}$ such that \begin{equation}
\bT^t\bF(\overline{\bx})\bD^{\frac{1}{2}}=\bU\bR\bE,\end{equation}
where $\bR\in \real^{n\times(m+p)}$ and
$$\bR_{ij}=\left\{\begin{array}{ll}
                          \sqrt{a_i}, & i=j\text{ and }i=1,\ldots,{m+p} \\
                          0, & \text{otherwise.}
                        \end{array}
\right.$$
According to the singular value decomposition theory,
we know that $\bU$ is the identity matrix. Then
\begin{eqnarray*}\bgra^2\Pi^d(\overline{\bvsig}) & = &
-\bF(\overline{\bx})^t(\bG(\overline{\bvsig}))^{-1}\bF(\overline{\bx})-\bD^{-1}\\
& = &
-(\bF(\overline{\bx})^t\bT)[\bT^t\bG(\overline{\bvsig})\bT]^{-1}(\bT^t\bF(\overline{\bx}))-\bD^{-1}\\
& = &
-\bD^{-\frac{1}{2}}\bE^t\bR^t\Diag\left(-\frac{1}{\lam_1},\ldots,-\frac{1}{\lam_n}\right)\bR\bE\bD^{-\frac{1}{2}}-\bD^{-1}\succeq
0. \end{eqnarray*} Multiplying by $\bD^{\frac{1}{2}}$ from the left and the right
\eb
\bD^{\frac{1}{2}}\bgra^2\Pi^d(\overline{\bvsig})\bD^{\frac{1}{2}}
=
-\bE^t\bR^t\Diag\left(-\frac{1}{\lam_1},\ldots,-\frac{1}{\lam_n}\right)\bR\bE
- {\bf I}_{(m+p)\times(m+p)}\succeq 0.
\ee
If we multiply the right side of the last equation by $\bE$ from the left and $\bE^t$ from the right,
 we have \begin{eqnarray*}0 & \preceq &
-\bR^t\Diag\left(-\frac{1}{\lam_1},\ldots,-\frac{1}{\lam_n}\right)\bR - {\bf I}_{(m+p)\times(m+p)}\\
& \preceq & \Diag\left(\frac{a_1}{\lam_1}-1,\ldots,\frac{a_{m+p}}{\lam_{m+p}}-1\right),\end{eqnarray*} thus $a_i \geq\lam_i,$
for every $i=1,\ldots,{m+p}.$ On the other hand \begin{eqnarray*}  \bT^t\bgra^2\Pi(\overline{\bx})\bT & = & \bT^t\bG(\overline{\bvsig})\bT+\bT^t\bF(\overline{\bx})\bD\bF(\overline{\bx})^t\bT\\
& = & \Diag(-\lam_1,\ldots,-\lam_n)+ \Diag(a_1,\ldots,a_{m+p},0,\ldots,0)\\
& = & \Diag(a_1-\lam_1,\ldots,a_{m+p}-\lam_{m+p},-\lam_{m+p+1},\ldots,-\lam_n).\end{eqnarray*}
Let $\bJ\in\real^{n\times n}$ be defined by $$J_{ij}=\left\{\begin{array}{ll}
                          1, & i=j\text{ and }i=1,\ldots,{m+p} \\
                          0, & \text{otherwise.}
                        \end{array}
\right.$$ Then we have \begin{equation}\label{LHESSPIL}
\bJ^t\bT^t\bgra^2\Pi(\overline{\bx})\bT\bJ=\Diag(a_1-\lam_1,\ldots,a_{m+p}-\lam_{m+p})\succeq 0.
\end{equation} Let ${\bf L}=\bT\bJ$, clearly $\Rank({\bf L})=m+p$ and ${\bf L}^t\bgra^2\Pi(\overline{\bx}){\bf L}\succeq 0$, this completes the proof. \hfill $\blacksquare$\\

In a similar way, we can prove the following lemma.

\begin{lemma}\label{MatrixQ} Suppose that $m+p>n$, $\overline{\bvsig}\in\calS_a^-$ is a stationary point $\Pi^d$ and $\bar{\bx}=(\bG(\overline{\bvsig}))^{-1}{\bf f}$ is a local minimizer of
$\Pi$. Then, there exists a matrix $\bQ\in\real^{(m+p)\times n}$ with
$\Rank(\bQ)=n$ such that
\begin{equation}\label{localminm1m2<n}
\bQ^t\bgra^2\Pi^d(\overline{\bvsig})\bQ\succeq 0.\end{equation}
\end{lemma}

Let the $m+p$ column vectors of ${\bf L}$ be respectively as ${\bf l}_1,\ldots,{\bf l}_{m+p}$ and the $n$ column vectors of $\bQ$ be respectively as ${\bf q}_1,\ldots,{\bf q}_n$. Clearly, ${\bf
l}_1,\ldots,{\bf l}_{m+p}$ are $m+p$ independent vectors and ${\bf  q}_1,\ldots,{\bf q}_n$ are $n$ independent vectors. Now the subspaces $\mathcal{X}_b$ and $\mathcal{S}_b$ are defined as follows:\begin{eqnarray} \mathcal{X}_b& = & \left\{\bx\in\real^n:\bx=\overline{\bx}+\sum_{i=1}^{m+p}\upsilon_i{\bf l}_i,\{\upsilon_i\}_{i=1}^{m+p}\subset \real\right\}, \\
\mathcal{S}_b& = & \left\{\bvsig\in\real^{m+p}: \bvsig=\overline{\bvsig}+\sum_{j=1}^{n}\vartheta_j{\bf q}_j,\{\vartheta_j\}_{j=1}^{n}\subset \real\right\}.
\end{eqnarray}

Now we are ready to present the Refined Triality Theory.

\begin{thm}[Triality Theory]\label{ThExtrC}
Let $\overline{\bvsig}$ be a stationary point of $\Pi^d$ and $\overline{\bx}=(\bG(\overline{\bvsig}))^{-1}{\bf f}$.
Assume that $ \det (\nabla^2 \Pi(\overline{\bx}) )  \neq 0$.
\begin{enumerate}
\item[(i)] If $\overline{\bvsig}\in \calS_a^+$, then $\overline{\bvsig}$ is the only global maximizer of $\Pi^d$ in $\calS_a^+$ and $\overline{\bx}$ is the only global minimizer of $\Pi$.
\item[(ii)] If $\overline{\bvsig}\in \calS_a^-$, then $\overline{\bvsig}$ is a local maximizer of $\Pi^d$ in $\calS_a^-$ if and only if $\overline{\bx}$ is a local maximizer of $\Pi$.
\item[(iii)] If $\overline{\bvsig}\in \calS_a^-$ and
\begin{enumerate}
\item[a)] if $n=m+p$, then $\overline{\bvsig}$ is a local minimizer of $\Pi^d$ if and only if $\overline{\bx}$ is a local minimizer of $\Pi$, i.e., there exists respectively
neighborhoods $\mathcal{X},\mathcal{S}\subset\real^n$ of $\overline{\bx}$ and
$\overline{\bvsig}$  such that
\begin{equation}
\Pi(\overline{\bx})=\min_{\bx\in \mathcal{X}}\Pi(\bx)=\min_{\bvsig\in
\mathcal{S}}\Pi^d(\bvsig)=\Pi^d(\overline{\bvsig});
\end{equation}
\item[b)]
if $m+p<n$ and $\overline{\bvsig}$ is a local minimizer of $\Pi^d$, then $\overline{\bx}$ is a saddle point of $\Pi$ and there exists respectively neighborhoods $\mathcal{X},\mathcal{S}\subset\real^n$  of $\overline{\bx}$ and $\overline{\bvsig}$, such that
\begin{equation}\label{WeakMinDP}
\Pi(\overline{\bx})=\min_{\bx\in \mathcal{X}\cap \mathcal{X}_b}\Pi(\bx)=\min_{\bvsig\in \mathcal{S}}\Pi^d(\bvsig)=\Pi^d(\overline{\bvsig});
\end{equation}
\item[c)] if $n<m+p$ and $\overline{\bx}$ is a local minimizer of $\Pi$, then $\overline{\bvsig}$ is a saddle point of $\Pi^d$ and there exists respectively neighborhoods $\mathcal{X},\mathcal{S}\subset\real^n$ of $\overline{\bx}$ and $\overline{\bvsig}$  such that \begin{equation} \Pi(\overline{\bx})=\min_{\bx\in \mathcal{X}}\Pi(\bx)=\min_{\bvsig\in \mathcal{S}\cap  \mathcal{S}_b}\Pi^d(\bvsig)=\Pi^d(\overline{\bvsig}).\end{equation}
\end{enumerate}
\end{enumerate}
\end{thm}

\textbf{Proof:}

(i) Since $\overline{\bvsig}\in \calS_a^+$, from Equation \eqref{HessianPid} it is not difficult to show that $\Pi^d$ is strictly concave in $\calS_a^+$ and $\Xi(\cdot,\overline{\bvsig})$ is strictly convex in $\real^n$ and therefore $\overline{\bvsig}$ must be the only global maximizer of $\Pi^d$ in $\calS_a^+$ and $\overline{\bx}$ is the only global minimizer of $\Xi(\cdot,\overline{\bvsig})$. By the definition of $\Xi$ given in Equation
\eqref{xireduced} and the convexity of $V$,  the Fenchel inequality leads to
$$\Xi(\bx,\bvsig)\leq \Pi(\bx),\ \forall (\bx,\bvsig)\in \real^n\times \calS_a.$$

Let us assume now that there exists a vector $\bx'\in \real^n\setminus\{\overline{\bx}\}$ such that $\Pi(\bx')\leq\Pi(\overline{\bx})$, then $$\Pi(\overline{\bx})\geq\Pi(\bx')\geq\Xi(\bx',\overline{\bvsig})>
\Xi(\overline{\bx},\overline{\bvsig})=\Pi(\overline{\bx}),$$ where the last equality comes from Equation \eqref{pixeqspidsig}. This contradiction proves that $\overline{\bx}$ must be the only global minimizer of $\Pi$.

(ii)  Notice first that using Equations \eqref{sigmaandlambdap1} and \eqref{sigmaandlambdap2} in Equation \eqref{HessianPi} we have  \begin{equation}\label{HessianPi2}\bgra^2\Pi(\overline{\bx})=\bG(\overline{\bvsig})+\bF(\overline{\bx})\bD\bF
(\overline{\bx})^t,\end{equation} where $F(\bx)$ and $D$ are defined in Equation \eqref{HessianPid}.
    If $\overline{\bvsig}$ is a local maximizer of $\Pi^d$ in $\calS_a^-$ we must have that $\bgra^2\Pi^d(\overline{\bvsig})\preceq 0$, from Equation \eqref{HessianPid} which  is equivalent to \begin{equation}\label{locmaxpid}\bD^{-1}+\bF(\overline{\bx})^t(\bG(\overline{\bvsig}))^{-1}\bF(\overline{\bx})\succeq 0.
    \end{equation}

    \begin{enumerate}
    \item[$\bullet$] If $m+p=n$ and $\bF$ is invertible, multiplying Equation \eqref{locmaxpid} by $(\bF(\overline{\bx})^t)^{-1}$ from the left and $(\bF(\overline{\bx}))^{-1}$ from the right,
        we have: \begin{equation}(\bF(\overline{\bx})^t)^{-1}\bD^{-1}(\bF(\overline{\bx}))^{-1}+(\bG(\overline{\bvsig}))^{-1}\succeq 0\end{equation} this is equivalent to $$(\bF(\overline{\bx})^t)^{-1}\bD^{-1}(\bF(\overline{\bx}))^{-1}\succeq -(\bG(\overline{\bvsig}))^{-1}\succ 0,$$ which in turn is equivalent to (Lemma \ref{G>UiffU1<G1} in the Appendix) $$-\bG(\overline{\bvsig})\succeq \bF(\overline{\bx})\bD\bF(\overline{\bx})^t\Longleftrightarrow \bgra^2\Pi(\overline{\bx})\preceq 0.$$
        By assumption $ \det ( \nabla^2 \Pi(\bar\bx) ) \neq 0$, then
         $\overline{\bx}$ is a local maximum of $\Pi$.
    \item[$\bullet$] If $m+p\neq n$ or $\bF$ is not invertible, then by Lemma \ref{SVDLEMMA}, there exists orthogonal matrices $\bE\in\real^{n\times n}$, $\bK\in \real^{(m+p)\times(m+p)}$ and a matrix $\bR\in \real^{n\times(m+p)}$ such that $$R_{ij}=\left\{\begin{array}{ll}
                          s_i, & i=j\text{ and }i=1,\ldots,r \\
                          0, & \text{otherwise}
                        \end{array}
        \right.$$ where $s_i>0$ for every $i$, $r=\Rank(\bF(\overline{\bx}))$
        and
        \begin{equation}\label{FDeqsERK}
        \bF(\overline{\bx})\bD^{\frac{1}{2}}=\bE\bR\bK .
        \end{equation}
         Using Equation \eqref{FDeqsERK},
         Equation \eqref{locmaxpid} can be rewritten as:
\[
\bD^{-1}+\bD^{-\frac{1}{2}}\bK^t\bR^t\bE^t(\bG(\overline{\bvsig}))^{-1}\bE\bR\bK\bD^{-\frac{1}{2}}\succeq 0
\]
after multiplying this equation by $\bK\bD^{\frac{1}{2}}$
from the left and $\bD^{\frac{1}{2}}\bK^t$ from the right, we have
$${\bf I}_{(m+p)\times(m+p)}+\bR^t(\bE^t\bG(\overline{\bvsig})\bE)^{-1}\bR\succeq 0.$$
This equation is equivalent to
$$-{\bf I}_{(m+p)\times(m+p)}-\bR^t(\bE^t\bG(\overline{\bvsig})\bE)^{-1}\bR\preceq 0.$$ By  Lemma \ref{WuMatLemma} in the Appendix, the last equation is equivalent to $$0\succeq \bE^t\bG(\overline{\bvsig})\bE+\bR\bR^t
=\bE^t\bG(\overline{\bvsig})\bE+\bR(\bK\bD^{-\frac{1}{2}}\bD\bD^{-\frac{1}{2}}\bK^t)\bR^t$$
 multiplying by $\bE$ from the left and $\bE^t$ from the right, we can obtain that
$$0\succeq \bG(\overline{\bvsig})+(\bE\bR\bK\bD^{-\frac{1}{2}})\bD(\bD^{-\frac{1}{2}}\bK^t\bR^t\bE^t)=
\bG(\overline{\bvsig})+\bF(\overline{\bx})\bD\bF(\overline{\bx})^t= \nabla^2 \Pi(\overline{\bx}).$$
By the assumption $ \det (\nabla^2 \Pi(\barbx))  \neq 0$,
$\overline{\bx}$ is a local maximum of $\Pi$.\\
\end{enumerate}
Notice that every step of the proof is  equivalent, so if $\overline{\bx}$ is a local maximum of $\Pi$ then $\overline{\bvsig}$ must be a local maximum of $\Pi^d$.

(iii) Let us consider the three cases:
\begin{enumerate}
\item[a)] $n=m+p$: if $\overline{\bvsig}$ is a local minimizer of $\Pi^d$ then
$$  \bgra^2\Pi^d(\overline{\bvsig}) =-\bF(\overline{\bx})^t(\bG(\overline{\bvsig}))^{-1}\bF(\overline{\bx})- \bD^{-1}\succeq 0 $$
$$   \Longleftrightarrow  -\bF(\overline{\bx})^t(\bG(\overline{\bvsig}))^{-1}\bF(\overline{\bx})\succeq \bD^{-1}.
$$
 This implies that
     $\Rank(\bF(\overline{\bx}))=n$.
      By multiplying the last inequality by $(\bF(\overline{\bx})^t)^{-1}$ from the left and by $(\bF(\overline{\bx}))^{-1}$ from the right, we have $$-(\bG(\overline{\bvsig}))^{-1}\succeq (\bF(\overline{\bx})^t)^{-1}\bD^{-1} (\bF(\overline{\bx}))^{-1}.$$
      By Lemma \ref{G>UiffU1<G1} this is equivalent to $$-\bG(\overline{\bvsig})\preceq
\bF(\overline{\bx})\bD\bF(\overline{\bx})^t \Longleftrightarrow \bgra^2\Pi(\overline{\bx})\succeq 0.$$ And since  $ \det (\nabla^2 \Pi(\barbx))  \neq 0$, $\overline{\bx}$ is a local minimizer of $\Pi$. In a similar way we can prove the converse.
\item[b)] From Equation \eqref{HessianPid} we know that
$$-\bF(\overline{\bx})^t(\bG(\overline{\bvsig}))^{-1}\bF(\overline{\bx})\succeq
\bD^{-1},$$
then $-\bF(\overline{\bx})^t(\bG(\overline{\bvsig}))^{-1}\bF(\overline{\bx})$ is a nonsingular matrix and $\Rank(\bF(\overline{\bx}))=m+p<n$.
We claim now that $\overline{\bx}$ is not a local minimizer of $\Pi$.
This is because that  if $\overline{\bx}$ is also a local minimizer, we would have  $$\bgra^2\Pi(\overline{\bx})=\bG(\overline{\bvsig})+\bF(\overline{\bx})\bD\bF(\overline{\bx})^t\succeq 0,$$ thus $$\bF(\overline{\bx})\bD\bF(\overline{\bx})^t\succeq -\bG(\overline{\bvsig}).$$
This implies that $$n=\Rank(-\bG(\overline{\bvsig}))=\Rank(\bF(\overline{\bx})\bD\bF(\overline{x})^t)=m+p,$$
which is a contradiction. Therefore, $\overline{\bx}$ is a saddle point of $\Pi$.\\

To prove Equation \eqref{WeakMinDP}, we let ${\bf L} $ be the matrix as given
in Lemma \ref{MatrixL} and $\{{\bf l}_i\}_{i=1}^{m+p}$ be the column vectors of ${\bf L}$. Define
$$\varphi(t_1,\ldots,t_{m+p}):=\Pi(\overline{\bx}+t_1{\bf l}_1+\ldots+t_{m+p}{\bf l}_{m+p}).$$
We need to show that $(0,\ldots,0)\in\real^{m+p}$ is a local
minimizer of the function $\varphi$.
Notice that
$$\bgra\varphi(0,\ldots,0)={\bf L}^t\bgra\Pi(\overline{\bx})=0$$ and
$$\bgra^2\varphi(0,\ldots,0)={\bf L}^t\bgra^2\Pi(\overline{\bx}){\bf L}\succeq
0,$$ which is a consequence of Lemma \ref{MatrixL}.
 Furthermore, from Equation \eqref{LHESSPIL} we have that $$\bgra^2\varphi(0,\ldots,0)=\Diag(a_1-\lam_1,\ldots,a_{m+p}-\lam_{m+p}),$$ and since $ \det (\nabla^2 \Pi(\overline{\bx}) )  \neq 0$ it can be proven that $a_i>\lam_i$ for every $i$. The proof is complete.
\item[c)] The proof is similar with  item b).\hfill $\blacksquare$
\end{enumerate}

\begin{remark}\label{rem-2}
  Theorem \ref{ThExtrC} shows that in order to solve the problem $(\calP) $ by means of the canonical duality theory, a necessary condition is that the problem $(\calP)$ should have a unique solution.
It was indicated in \cite{ruan-gao-pe} that if the nonconvex minimization problem has more than one global minimizer, it could be NP-hard.
In order to solve this type of problems, the perturbation methods should be used.
\end{remark}

\begin{remark} The triality theory states precisely that if $\bvsig$ is a global maximizer of $\Pi^d$ on a certain set, then $\bx$ is a global minimizer for $\Pi$. This is known from the general result by Gao and Strang in \cite{gao-gs1}. If $\bvsig$ is a local maximizer for $\Pi^d$ then $\bx$ is also a local maximizer for $\Pi$. This is the so-called double-max duality statement. If $\bvsig$ is a local minimizer for $\Pi^d$, then $\bx$ is also a local minimizer for $\Pi$ in certain directions. This is so-called double-min duality in the standard triality form proposed in \cite{DGAO1}. The triality theory was first discovered in nonconvex mechanics \cite{gao-amr}. It was realized in 2003 that the double-min duality holds under certain additional condition, which was left as an open problem  (see \cite{Gao-amma03,gao-opt03}).
   Recently, this open problem is solved for quartic polynomial optimization problem \cite{WU}. This result is now generalized to the general nonconvex problem $(\calP)$. Part (iii)  of Theorem \ref{ThExtrC} shows  that if $m+p=n$, then $\bvsig$  is a local minimizer if and only if $\bx$ is also a local minimizer. In other cases either $\bx$ is a saddle point of $\Pi$ or $\bvsig$ is a saddle point of $\Pi^d$.
\end{remark}
\begin{remark} 
 The canonical duality-triality theory has been challenged recently by   C. {\em $Z\check{a}$}linescu and his co-workers
R. Strugariu,  M. D. Voisei
  in several papers  (see \cite{v-z}).
   By list some simple ``counterexamples", they claimed that this  theory is false.
   Unfortunately, most of these counterexamples are  not new, which were first discovered by Gao  in 2003 \cite{Gao-amma03,gao-opt03}.
  However,  \cite{Gao-amma03,gao-opt03}   never been cited in their papers.
   Some  of their ``counterexamples" are  fundamentally wrong, i.e. they oppositely choose 
   linear functions as the stored energy and nonlinear  functions as external energy (see \cite{vz2}). 
   These conceptual mistakes show a big gap between mathematics and mechanics.
   \end{remark} 
   
\section{Numerical Examples}

In the following examples, $m=p=1$ and $\beta_1=1$. The graphs provided and the numerical results were obtained using Maxima \cite{MAXIMA}.

\subsection{  One stationary point in $\calSa^+$}

First, we consider the case that the primal function has a unique solution. We let $ \alpha_1=\theta_1=1$ and
$$\bA=\left[\begin{array}{cc}
                                          1 & 0 \\
                                          0 & -1
                                        \end{array}
\right], \bB_1=\left[\begin{array}{cc}
                     1 & 0 \\
                     0 & 2
                   \end{array}
\right], \bC_1=\left[\begin{array}{cc}
                     1 & 0 \\
                     0 & 1
                   \end{array}
\right], {\bf f}=\left[\begin{array}{c}
                   1 \\
                   1
                 \end{array}
\right] .$$
Clearly, the function $\Pi:\real^2\rightarrow \real$ is given by
$$\Pi(x,y) =\exp\left(\frac{1}{2}(x^2+2y^2)-1\right)+\frac{1}{2}\left(\frac{1}{2}(x^2+y^2)-1\right)^2+\frac{1}{2}(x^2-y^2)-x-y,$$
and the dual function has the form of
$$\Pi^d(\tau,\sig) =-\frac{1}{2}\left(\frac{1}{1+\tau+\sig}+\frac{1}{2\tau+\sig-1}\right)-\tau\cdot
\ln(\tau)-\frac{1}{2} \sig^2-\sig.$$

\begin{figure}[h]
  \centering
  \subfloat[Contour Levels of $\Pi$.]{\label{gex1a}\includegraphics[width=0.4\textwidth]{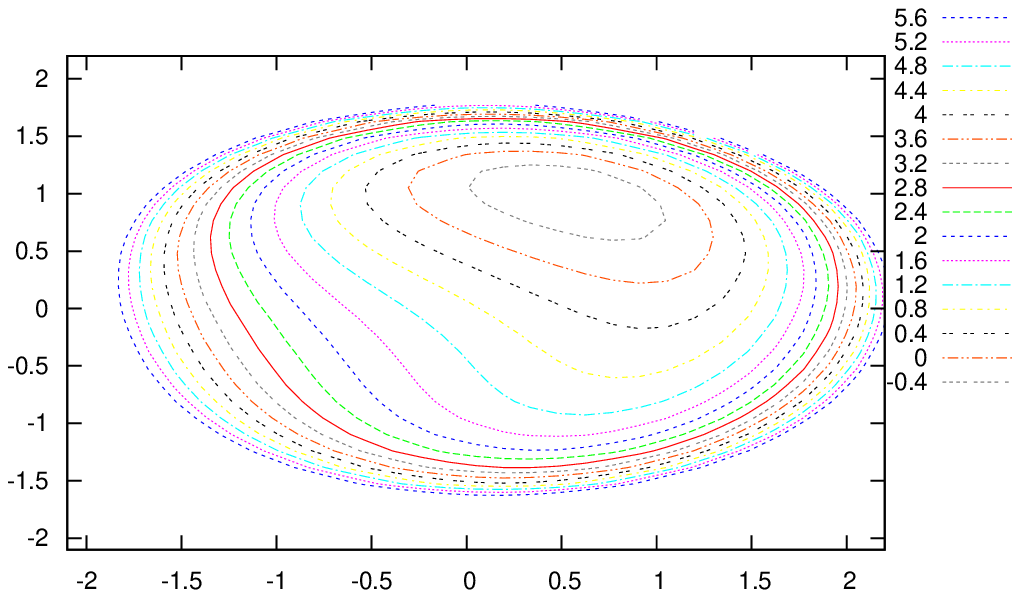}}\quad
  \subfloat[Graph of $\Pi$.]
  {\label{gex1b}\includegraphics[width=0.4\textwidth]{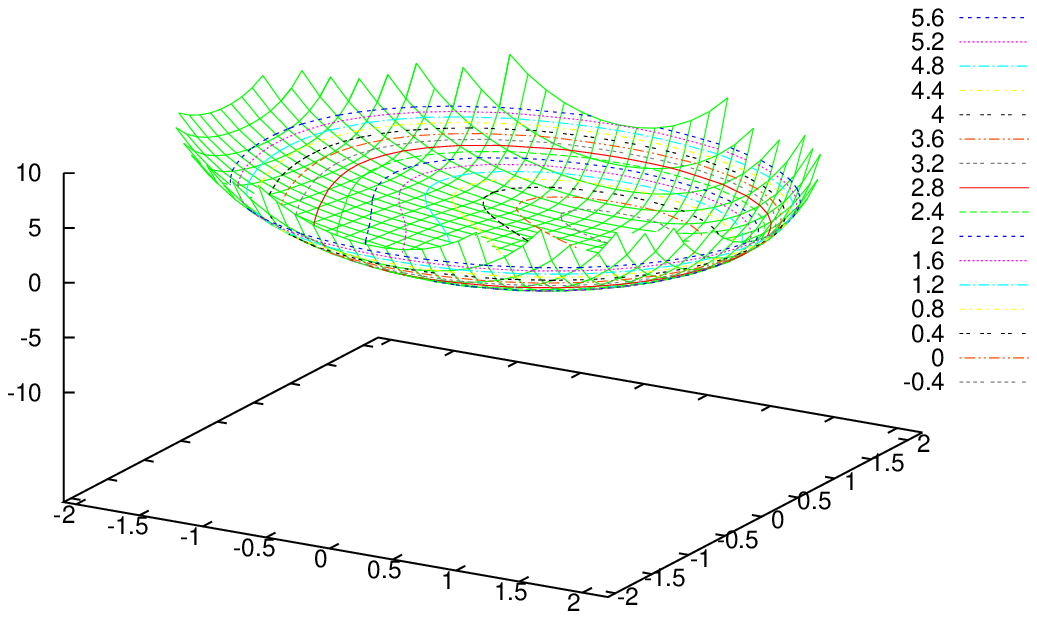}}
  \caption{$\Pi$ function of Example 1} \label{gex1}
\end{figure}

It can be shown that $\Pi^d$ has only one critical point in
$\calS_a^+$ and it is given (approximately) by
$$\overline{\bvsig}=(1.171057661103504,-0.34599084656216). $$
By the triality theory, the vector
 \[
 \overline{\bx}=\bG(\overline{\bvsig})^{-1}{\bf f} =(0.54792514555217,1.003890602479819)
 \]
 is the only global minimizer of the primal problem.

\begin{figure}
  \centering
  \subfloat[Contour Levels of $\Pi^d$.]{\label{gex1c}\includegraphics[width=0.4\textwidth]{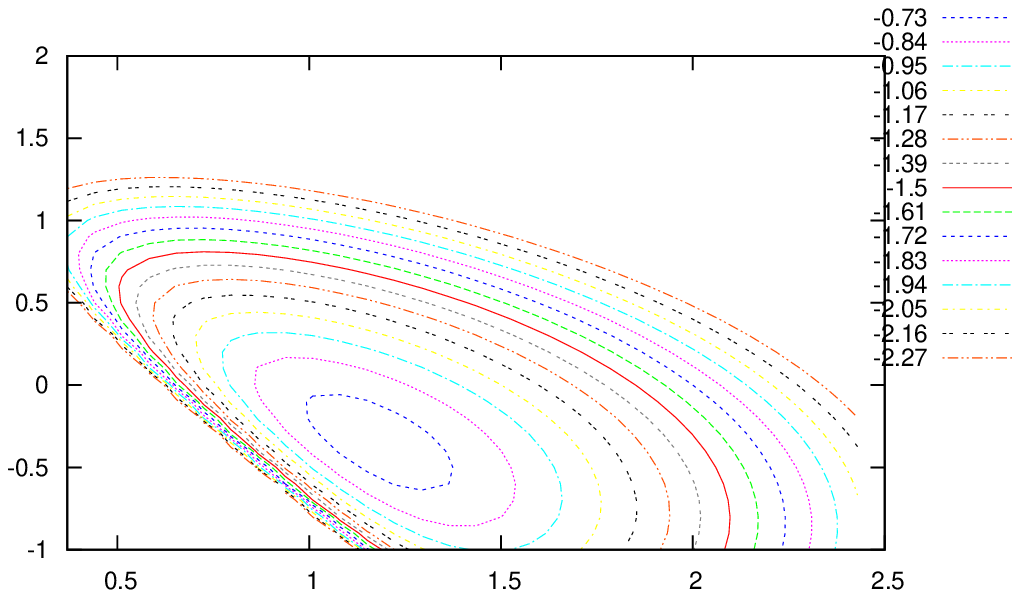}}\quad
  \subfloat[Graph of $\Pi^d$.]
  {\label{gex1d}\includegraphics[width=0.4\textwidth]{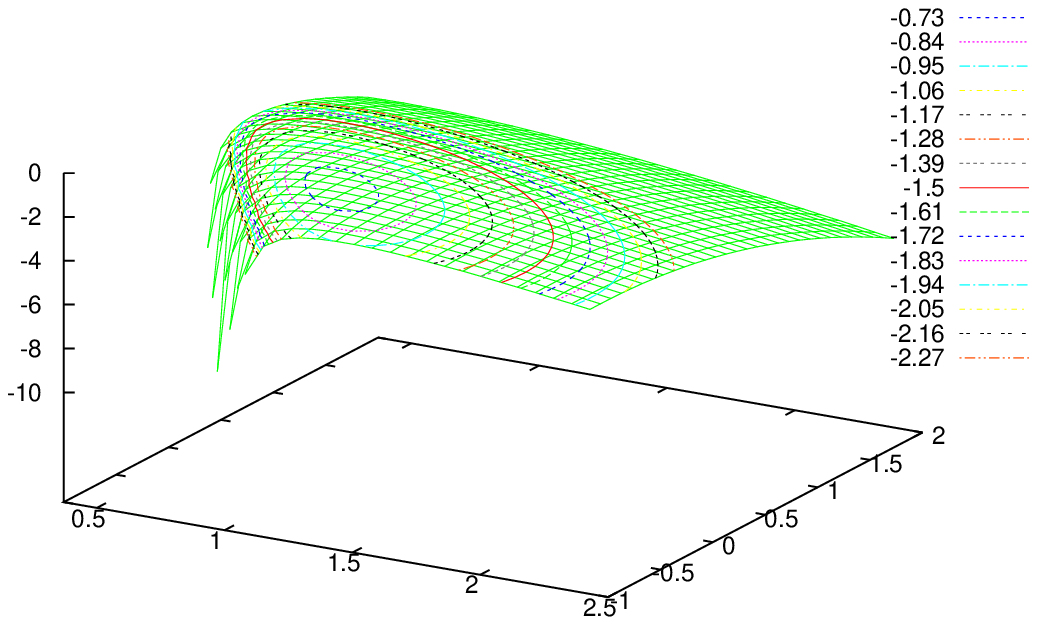}}
  \caption{$\Pi^d$ function of Example 1} \label{gex1D}
\end{figure}

\subsection{One stationary point in $\calSa^+$ and one   in $\calSa^-$}

\begin{figure}
  \centering
  \subfloat[Contour Levels of $\Pi$.]{\label{gex3a}\includegraphics[width=0.4\textwidth]{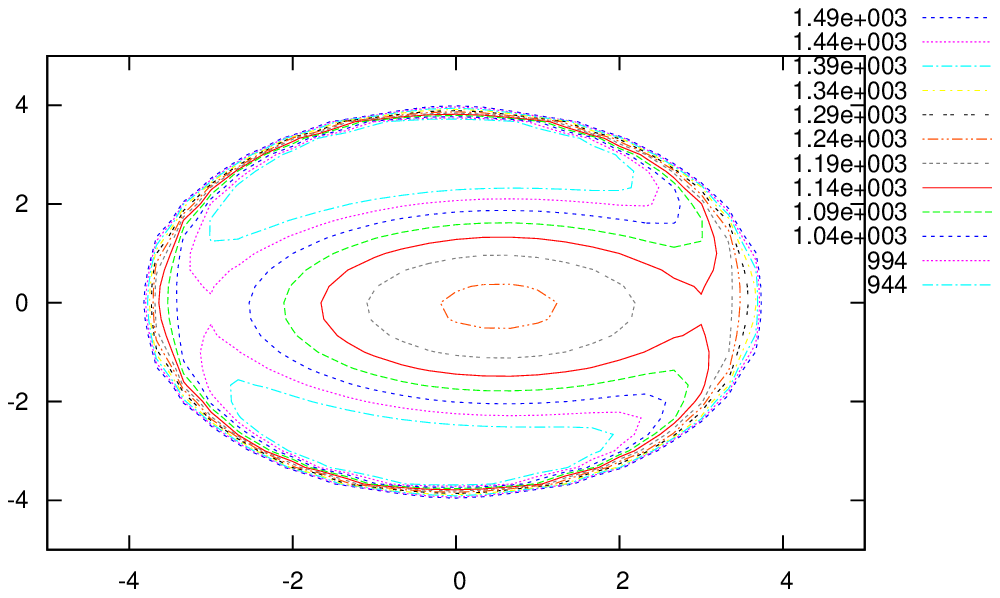}}\quad
  \subfloat[Graph of $\Pi$.]
  {\label{gex3b}\includegraphics[width=0.4\textwidth]{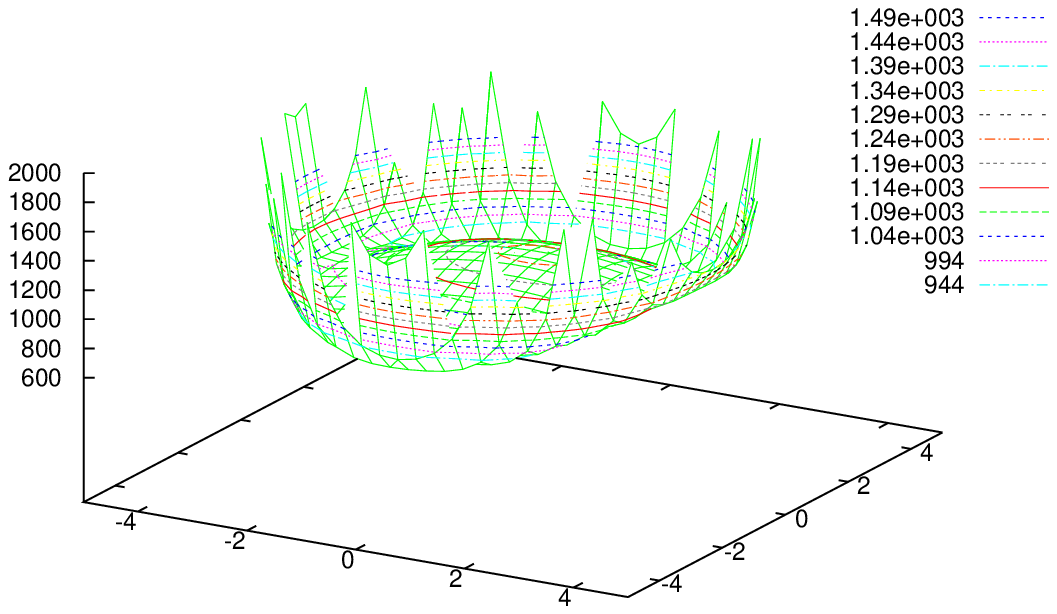}}
  \caption{Primal function $\Pi$  in  Example 2} \label{gex3}
\end{figure}

\begin{figure}
  \centering
  \subfloat[Contour Levels of $\Pi^d$ in $\calS_a^-$.]{\includegraphics[width=0.4\textwidth]{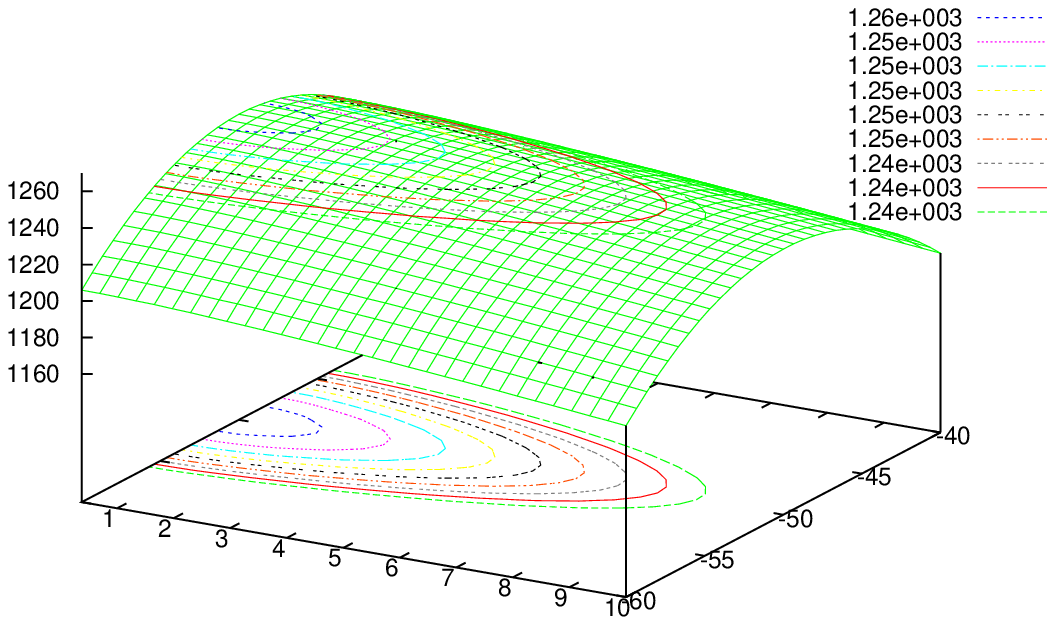}}\quad
  \subfloat[Graph of $\Pi^d$ in $\calS_a^-$.]
  {\includegraphics[width=0.4\textwidth]{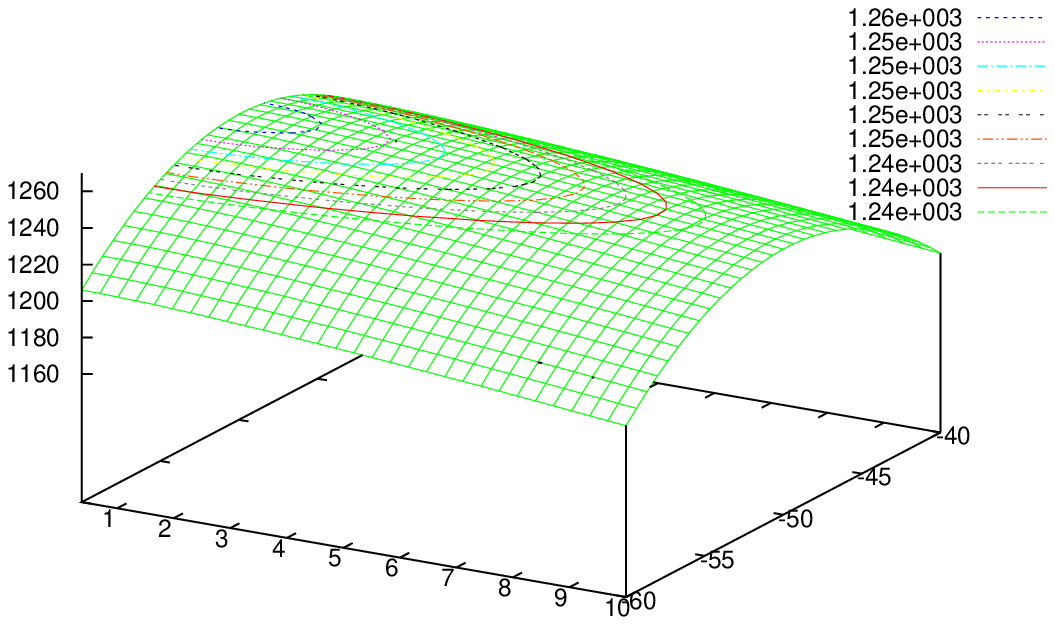}}
  \caption{$\Pi^d$ function in $\calS_a^-$ of Example 2}\label{gex3locmax}
\end{figure}

Let  $ \alpha_1=1, \;\; \theta_1=50$, and
$$\bA=\left[\begin{array}{cc}
                                          1 & 0 \\
                                          0 & -16
                                        \end{array}
\right], \bB_1=\left[\begin{array}{cc}
                     1 & 0 \\
                     0 & 1
                   \end{array}
\right], \bC_1=\left[\begin{array}{cc}
                     1 & 0 \\
                     0 & 2
                   \end{array}
\right], {\bf f}=\left[\begin{array}{c}
                   -25 \\
                   9
                 \end{array}
\right] .$$
The primal function $\Pi:\real^2\rightarrow \real$ is then given  by
$$\Pi(x,y) =\exp\left(\frac{1}{2}(x^2+y^2)-1\right)+\frac{1}{2}\left(\frac{1}{2}(x^2+2y^2)-50\right)^2
+\frac{1}{2}(x^2-16y^2)+25x-9y $$
and its canonical dual  is
$$\Pi^d(\tau,\sig) =-\frac{1}{2}\left(\frac{81}{-16+\tau+2\sig}+\frac{625}{1+\tau+\sig}\right)-\tau\cdot
\ln(\tau)-\frac{1}{2} \sig^2-50\sig, $$
which has  two critical points:
\[
\overline{\bvsig_1}=(96.61711963278241,-38.94928057661689) \in \calS_a^+,
\]
\[
 \overline{\bvsig_2}=(0.42157060067968,-49.86072154366873) \in \calS^-_a.
 \]
 Therefore, by the triality theory, the associated vector
 \[
 \overline{\bx_1} = \bG(\overline{\bvsig_1})^{-1}{\bf f}=(-0.42612784793499,3.310578038951848)
 \]
 is the only global minimizer of $\Pi(\bx)$ and $$\overline{\bx_2}=(0.51611144112381,-0.078057328303129)$$ is a local maximizer (see Figure \ref{gex3}) since $\overline{\bvsig_2}$ is a local maximum of $\Pi^d$ in $\calS_a^-$ (see Figure \ref{gex3locmax}).

\begin{figure}
  \centering
  \subfloat[Contour Levels of $\Pi^d$ in $\calS_a^+$.]{\label{gex3c}\includegraphics[width=0.4\textwidth]{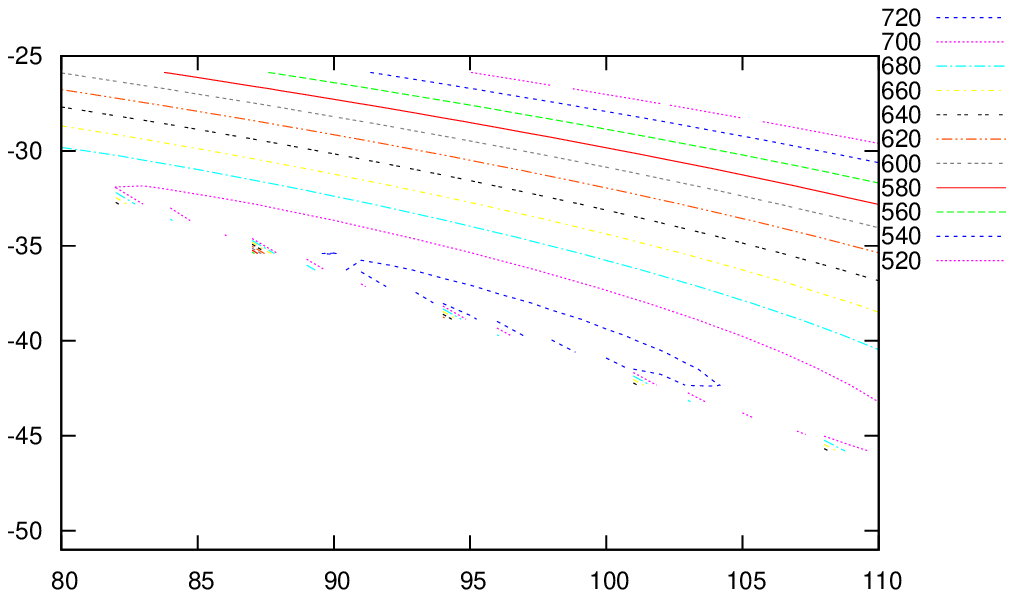}}\quad
  \subfloat[Graph of $\Pi^d$ in $\calS_a^+$.]
  {\label{gex3d}\includegraphics[width=0.4\textwidth]{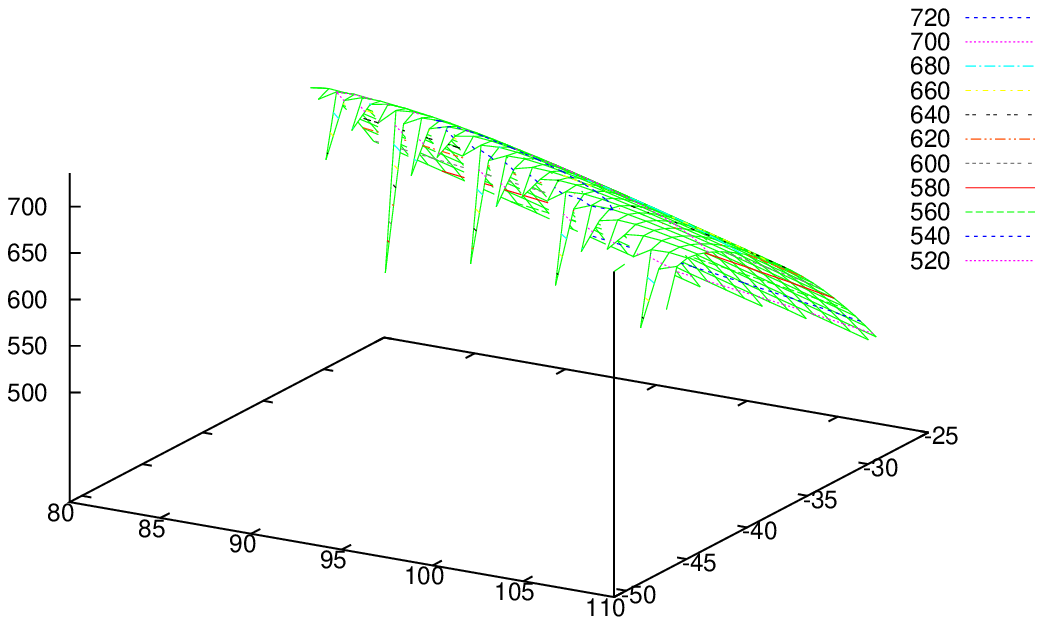}}
  \caption{$\Pi^d$ function in $\calS_a^+$ of Example 2} \label{gex3D}
\end{figure}

\subsection{One  stationary point in $\calSa^+$ and two  in $\calSa^-$}

In order to illustrate the triality theory, we let $ \alpha_1=\theta_1=2$, and
$$\bA=\left[\begin{array}{cc}
                                          -16 & 0 \\
                                          0 & -4
                                        \end{array}
\right], \bB_1=\left[\begin{array}{cc}
                     1 & 0 \\
                     0 & 0
                   \end{array}
\right], \bC_1=\left[\begin{array}{cc}
                     0 & 0 \\
                     0 & 1
                   \end{array}
\right], {\bf f}=\left[\begin{array}{c}
                   2 \\
                   2
                 \end{array}
\right] .$$
Accordingly,  we have
$$\Pi(x,y) =\exp\left(\frac{1}{2}x^2-2\right)+\frac{1}{2}\left(\frac{1}{2}y^2-2\right)^2+\frac{1}{2}(-16x^2-4y^2)-2x-2y,$$
$$\Pi^d(\tau,\sig) =-\frac{1}{2}\left(\frac{4}{\sig-4}+\frac{4}{\tau-16}\right)-\tau\cdot
\ln(\tau)-\tau-\frac{1}{2} \sig^2-2\sig.$$

\begin{figure}
  \centering
  \subfloat[Contour Levels of $\Pi$.]{\label{gex4a}\includegraphics[width=0.4\textwidth]{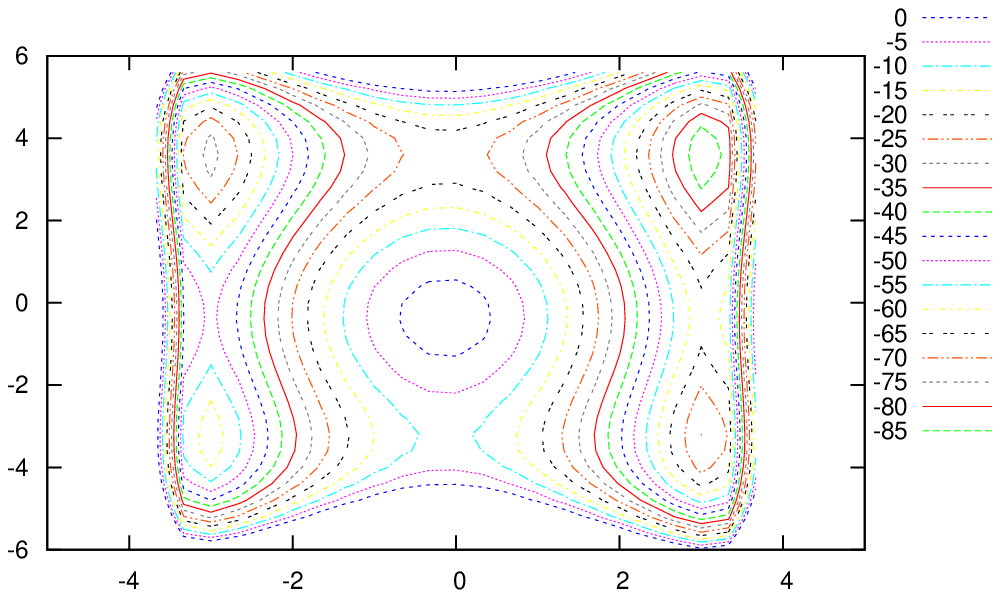}}\quad
  \subfloat[Graph of $\Pi$.]
  {\label{gex4b}\includegraphics[width=0.4\textwidth]{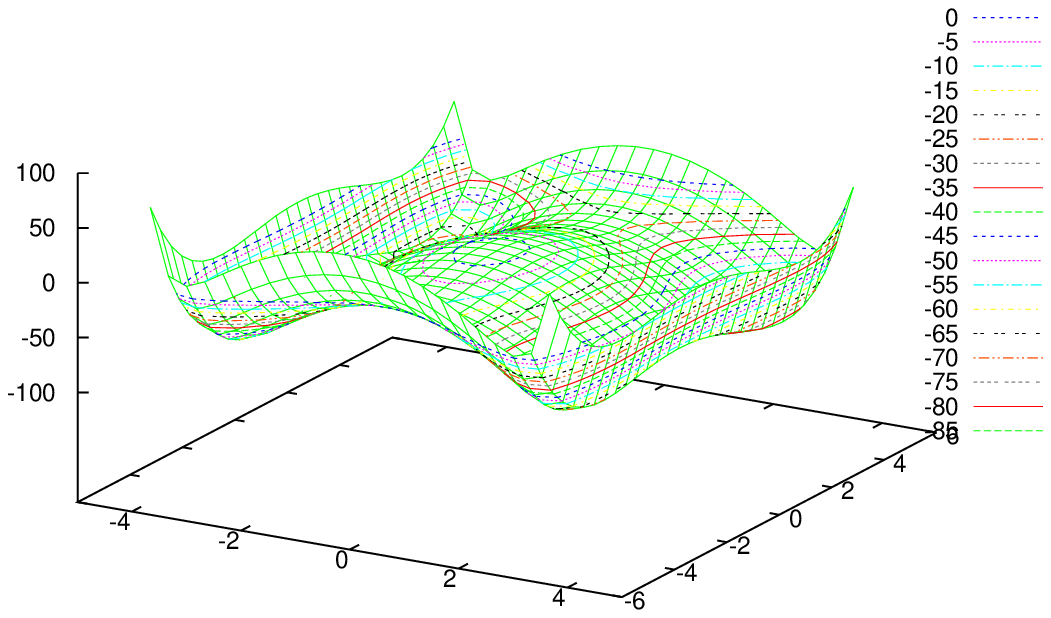}}
  \caption{$\Pi$ function of Example 3} \label{gex4}
\end{figure}

\begin{figure}
  \centering
  \subfloat[Contour Levels of $\Pi^d$ in $\calS_a^+$.]{\label{gex4c}\includegraphics[width=0.4\textwidth]{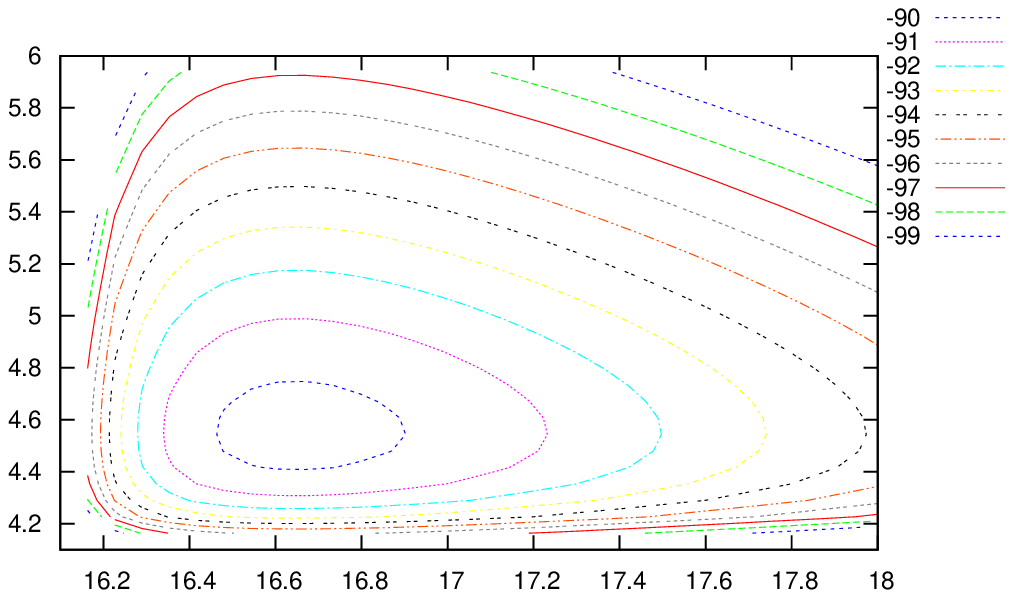}}\quad
  \subfloat[Graph of $\Pi^d$ in $\calS_a^+$.]
  {\label{gex4d}\includegraphics[width=0.4\textwidth]{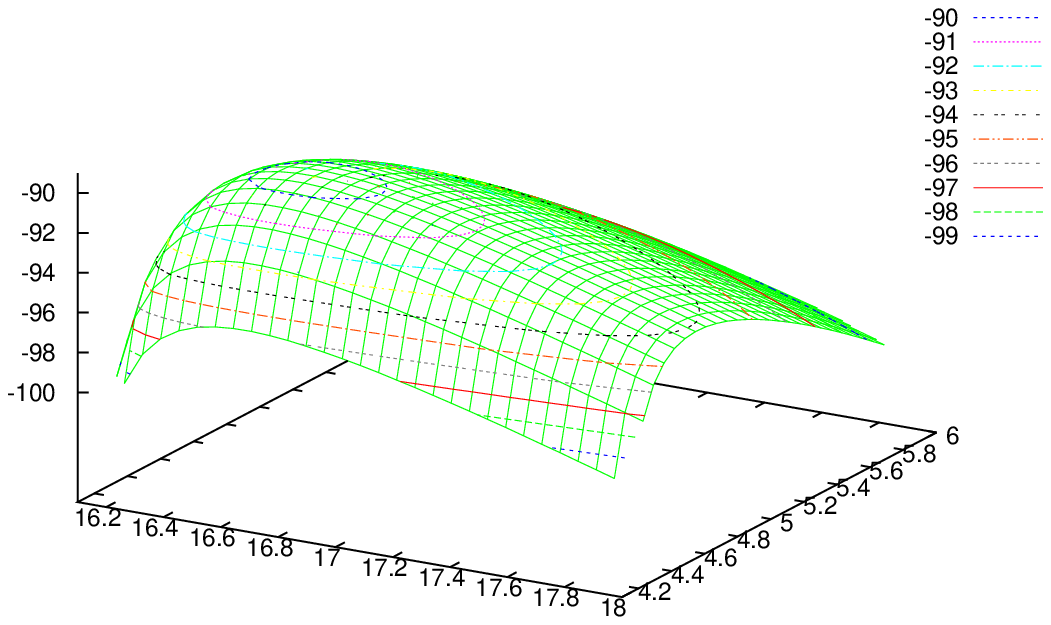}}
  \caption{$\Pi^d$ function in $\calS_a^+$ of Example 3}\label{gex4D}
\end{figure}

\begin{figure}
  \centering
  \subfloat[Contour Levels of $\Pi^d$ in $\calS_a^-$.]{\includegraphics[width=0.4\textwidth]{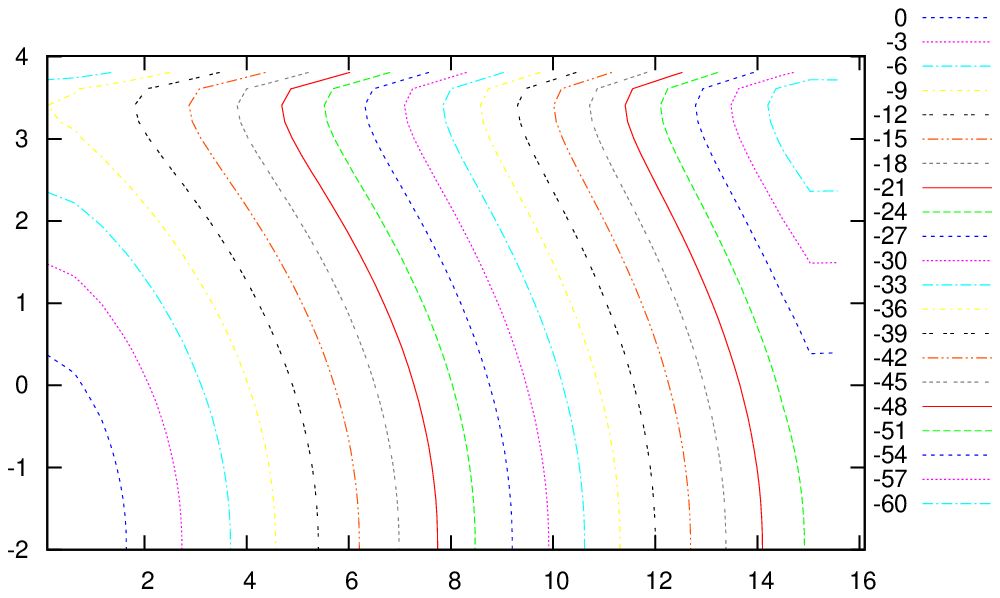}}\quad
  \subfloat[Graph of $\Pi^d$ in $\calS_a^-$.]
  {\includegraphics[width=0.4\textwidth]{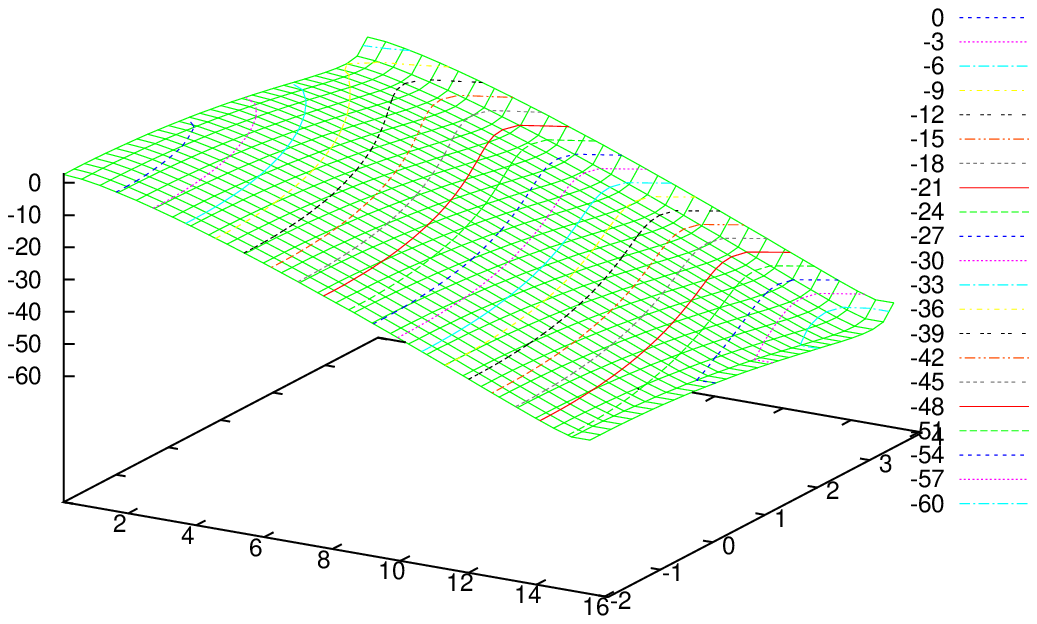}}
  \caption{$\Pi^d$ function in $\calS_a^-$ of Example 3}\label{gex4locmaxmin}
\end{figure}

In this case, $\Pi^d$ has in total six critical points but only one
\[
\overline{\bvsig_1}=(16.64468576727409,4.552474610531074) \in  \calS_a^+,
\]
(see Figure \ref{gex4D}) and two
\[
\overline{\bvsig_2}=(0.13641513779858,-1.943380912562619) \in \calS^-_a,
\]
\[
 \overline{\bvsig_3}=(15.34981976568548,3.390906302031545) \in \calS^-_a.
 \]
 From Figures
  \ref{gex4locmaxmin}  we can see that $\overline{\bvsig_2}$ is
a local maximizer   and $\overline{\bvsig_3}$
is a local minimizer  of $\Pi^d$. Therefore,  by the triality theory, we know that
\[
 \overline{\bx_1}=\bG(\overline{\bvsig_1})^{-1}{\bf f}=(3.102286573591542,3.620075858467906)
 \]
 is the only global minimizer;
\[
\overline{\bx_2}=(-0.12607490787063,-0.33650880356205)
\]
is a local maximizer  and
\[
\overline{\bx_3}=(-3.076070133243102,-3.283567054905852)
\]
is a local minimizer of $\Pi(\bx)$ (see  Figure
\ref{gex4}).

\subsection{Non-unique global minima}

In the case that no stationary point can be found in $\calS^+_a$, the primal problem could have more than one global minima. To see this, we let $ {\bf f}\equiv 0, \; \alpha_1=\theta_1=2$, and
 $$\bA\equiv 0, \bB_1=\left[\begin{array}{cc}
                     1 & 0 \\
                     0 & 0
                   \end{array}
\right], \bC_1=\left[\begin{array}{cc}
                     0 & 0 \\
                     0 & 1
                   \end{array}
\right] .$$
 In this case, the primal function
$$\Pi(x,y) =\exp\left(\frac{1}{2}x^2-2\right)+\frac{1}{2}\left(\frac{1}{2}y^2-2\right)^2 $$
has 2 global minimums at $(0,-2)$, $(0,2)$ and a local maximum at $(0,0)$. While the dual function $$\Pi^d(\tau,\sig) =-\tau\ln\tau -\tau -\frac{1}{2}\sig^2-2\sig $$ does not have a stationary point in $\calS_a^+$. There is however a critical point in the boundary of $\calSa^+$, namely, $\overline{\bvsig}=(\exp(-2),0)$. By defining $\overline{\bx}=\bG(\overline{\bvsig})^{-1}{\bf f}$, we have that $\overline{\bx}=(0,0)$.

In order to find a global minimum of $\Pi$, we need to introduce the following perturbations:
$$\bA_n=\left[\begin{array}{cc}
                                          -\frac{16}{n} & 0 \\
                                          0 & -\frac{4}{n}
                                        \end{array}
\right] \text{ and } {\bf f_n}=\left[\begin{array}{c}
                   \frac{2}{n} \\
                   \frac{2}{n}
                 \end{array}
\right],\text{ for every }n\in I\!\!N.$$
Then, the associated
 primal and dual functions are
\[
\Pi_n(x,y) =\exp\left(\frac{1}{2}x^2-2\right)+\frac{1}{2}\left(\frac{1}{2}y^2-2\right)^2
+\frac{1}{2}\left(-\frac{16}{n}x^2-\frac{4}{n}y^2\right)-\frac{2}{n}x-\frac{2}{n}y,
\]
\[
\Pi^d_n(\tau,\sig) = -\frac{1}{2}\left(\frac{4}{n^2\left(\tau-\frac{16}{n}\right)}+\frac{4}{n^2\left(\sig-\frac{4}{n}\right)}\right)-\tau\ln\tau+\tau-\frac{1}{2}\sig^2 -2\tau-2\sig.
\]
Notice that if $n=1$ we are in the case presented in Example 3.
Let us show that for sufficiently large values of $n$ we can find a stationary point for $\Pi^d_n$ in $\calS_a^+$, namely $\overline{\bvsig_n}$. Furthermore, by defining $\overline{\bx}_n=\bG(\overline{\bvsig}_n)^{-1}{\bf f_n}$ we will have a convergent sequence.

Let us calculate the gradient of $\Pi^d_n$:
\[
\bgra\Pi^d_n(\tau,\sig) =\left[\begin{array}{c}
 -2-\ln\tau+\frac{2}{(n\tau-16)^2} \\
  -\sig-2+\frac{2}{(n\sig-4)^2}
  \end{array}\right].
  \]

Let $h(\tau) =-2-\ln\tau+\frac{2}{(n\tau-16)^2}$ and
$g(\sig) = -\sig-2+\frac{2}{(n\sig-4)^2}$.
 It is not difficult to show that there exists a sufficiently large
 $N\in I\!\!N$, such that if $n>N$, the following are true:

\begin{enumerate}
\item[a)] $\displaystyle n\cdot \exp\left(-2+\frac{1}{n}\right)-16 \text{ and } n\cdot \exp\left(-2\right)-16$ are positive numbers.
\item[b)] $\displaystyle h\left(\exp\left(-2+\frac{1}{n}\right)\right)= \frac{2}{(n\cdot \exp(-2+\frac{1}{n})-16)^2} -\frac{1}{n}\\
     < 0 < h(\exp(-2))= \frac{2}{(n\cdot \exp(-2)-16)^2}$.
\item[c)] $\displaystyle g\left(\frac{5.1}{n}\right)\approx-\frac{5.1}{n}-0.34710743801< 0 < g\left(\frac{4.9}{n}\right)\approx0.46913580247-\frac{4.9}{n}.$
\end{enumerate}

Based on these results, we know that  for every $n>N$, $\bgra\Pi^d_n$ has a stationary point $\overline{\bvsig}_n=(\overline{\tau}_n,\overline{\sig}_n)\in$ $[\exp(-2),\exp(-2+\frac{1}{n})]\times\left[\frac{4.9}{n},\frac{5.1}{n}\right]$. Moreover, by
 the fact that $g(\overline{\sig}_n)=0$, it is easy to obtain
  $\displaystyle\lim_{n\rightarrow +\infty}n\cdot\overline{\sig}_n=5$.\\

Notice also that $$\bG(\overline{\bvsig}_n) =\left[\begin{array}{cc}
                                              \displaystyle\overline{\tau}_n-\frac{16}{n} & 0 \\
                                              0 & \overline{\sig}_n-\frac{4}{n}
                                            \end{array}\right] $$
                                            is positive definite.
                                            Therefore, the perturbed solution
                                            can be obtained as
\[
\overline{\bx}_n=\bG(\overline{\bvsig}_n)^{-1}{\bf f_n}=\left[\begin{array}{c}
  \frac{2}{n\cdot\overline{\tau}_n-16} \\
           \frac{2}{n\cdot\overline{\sig}_n-4}
  \end{array}
\right].
 \]
 Since $\overline{\tau}_n\in [\exp(-2),\exp(-2+\frac{1}{n})]$ then $\displaystyle\lim_{n\rightarrow +\infty}\tau_n=\exp(-2)$.
 By the fact that  $\displaystyle\lim_{n\rightarrow +\infty}n\cdot\overline{\sig}_n=5$,
  we have $$\lim_{n\rightarrow +\infty}\overline{\bx}_n=\left[\begin{array}{c}
                     0 \\
  2
   \end{array}
\right],$$ which is a solution of $\Pi$.

Canonical perturbation method was originally introduced in
\cite{r-g-j} for solving nonconvex polynomial minimization problems.
This method has been used successfully in integer programming and  network communication
(see \cite{gao-ruan-pardalos,wang-etal}).


\section{Future Research}

Some open questions that will be studied  in the future are the following:

\begin{enumerate}
\item[$\bullet$] As stated in Remark 1, in order to use the canonical dual transformation a necessary condition is that $(\calP)$ has a unique solution. Is this also a sufficient condition? In other words, giving $(\calP)$ such that it has a unique solution, can we find a stationary point of $\Pi^d$ in $\calSa^+$?
\item[$\bullet$] Example 4 shows an interesting perturbation method that allows us to solve a problem when the necessary condition of Remark 1 is not satisfied. Can we generalize this method and develop an algorithm?
\end{enumerate}

\section{Appendix: Some Lemmas in Matrix Analysis}

The following results are needed in the proofs of Section 2.

\begin{lemma}\label{SVDLEMMA} (Singular value decomposition \cite{HORN}) For any given matrix $\bM\subset\real^{m\times n}$ with $\Rank(\bM)=r$, there exists $\bU\subset\real^{m\times m}$, $\bR\subset\real^{m\times n}$ and $\bE\subset\real^{n\times n}$ such that
$$\bM = \bU\bR\bE;$$ where $\bU$ and $\bE$ are orthogonal matrices, and $$R_{ij}=\left\{\begin{array}{ll}
                                                              s_i, & i=j,\ i=1,\ldots,r \\
                                                              0, & i\neq j,
                                                            \end{array}\right.$$ where $s_i>0$ for every $i=1,\ldots,r$.
\end{lemma}

\begin{lemma}\label{G>UiffU1<G1} \cite{HORN} If $\bG$ and $\bU$ are positive definite matrices in $\real^{n\times n}$, then $\bG\succeq \bU$ if and only if $\bU^{-1}\succeq \bG^{-1}$.
\end{lemma}

\begin{lemma}\label{WuMatLemma} \cite{WU} Suppose $\bP$, $\bU$ and $\bD$ are three matrices in $\real^{n\times n}$ such that $$\bD=\left[\begin{array}{cc}
       \bD_{11} & 0_{m\times(n-m)} \\
    0_{(n-m)\times n} & 0_{(n-m)\times(n-m)}
 \end{array}
\right],$$ where $\bD_{11}\in \real^{m\times m}$ is nonsingular and
$$\bP=\left[\begin{array}{cc}
                                                                              \bP_{11} & \bP_{12} \\
                                                                              \bP_{21} & \bP_{22}
                                                                            \end{array}
\right]\prec 0,\ \bU=\left[\begin{array}{cc}
                           \bU_{11} & 0_{m\times(n-m)} \\
                           0_{(n-m)\times m} & \bU_{22}
                         \end{array}
\right] \succ 0,$$ $\bP_{ij}$ and $\bU_{ii}$ are appropriate
dimensional matrices for $i,j=1,2$. Then,
\begin{equation}\label{WuMatLemmaeq} \bP+\bD\bU\bD^t\preceq 0
\Longleftrightarrow -\bD^t\bP^{-1}\bD-\bU^{-1}\preceq
0.\end{equation}
\end{lemma}

$\;$\vspace{.5cm}\\

\noindent {\bf Acknowledgements}.\\
This  research is supported by US Air Force
Office of Scientific Research under the grant AFOSR FA9550-10-1-0487.
Comments and suggestions from editor and reviewers are sincerely acknowledged.


\begin{thebibliography}{99}\rm


\bibitem{asp-god}
Aspnes, J.;  Goldberg, D.; Yang, Y. R.
\newblock{\em On the computational complexity of sensor network localization.}
\newblock Lecture Notes in Computer Science (3121), Springer-Verlag, pp. 32-44 (2004).

\bibitem{cai-gao-qin}
Cai, K; Gao, D. Y. and Qin, Q. H.
\newblock {\em Post-buckling solutions of hyper-elastic beam by canonical dual finite element method.}
\newblock To appear in {\em Mathematics and Mechanics of Solids}.

\bibitem{DESOER}
Desoer, C. A.; Whalen, B. H.
\newblock {\em A Note on Pseudoinverses.}
\newblock Journal of the Society for Industrial and Applied
Mathematics, Vol. 11, No 2, pp. 442-447 (1963).


\bibitem{FENGLIN}
Feng, J. M.; Lin, G. X.; Sheu, R. L.; Xia, Y.
\newblock {\em Duality and solutions for quadratic programming over
single non-homogeneous quadratic constraint.}
\newblock Journal of Global Optimization, Published online: 17 Nov 2010.


\bibitem{gao-amr}
Gao, D. Y.
\newblock  {\em  Dual extremum principles in finite deformation theory     with applications to post-buckling analysis of extended nonlinear beam theory.}
\newblock Applied Mechanics Reviews,  50 (11), S64-S71 (1997).

\bibitem{gao-jem}
Gao, D. Y.
\newblock {\em Complementary finite element method for finite deformation nonsmooth mechanics.}
\newblock J. Eng. Math., 30,  pp. 339-353 (1996).


\bibitem{gao-ima}
Gao, D. Y.
\newblock {\em Duality, triality and complementary extremum principles in nonconvex parametric variational problems with applications.}
\newblock IMA J. Appl. Math., \textbf{61}, pp. 199-235 (1998).


\bibitem{gao-mecc}
Gao, D. Y.
\newblock {\em General Analytic Solutions and Complementary Variational Principles for Large Deformation Nonsmooth Mechanics.}
\newblock Meccanica, 34, pp. 169-198 (1999).

\bibitem{DGAO1}
Gao, D. Y.
\newblock {\em Duality Principles in nonconvex systems. Theory
Methods and Applications.}
\newblock  Kluwer Academic Publishers, Dordrecht/Boston (2000).


\bibitem{gao-jogo00}
Gao, D. Y.
\newblock {\em Canonical dual transformation method and generalized triality theory in nonsmooth global optimization.}
\newblock J. Glob. Optim. 17(1/4), pp. 127-160 (2000).


\bibitem{Gao-amma03}
Gao, D. Y.
\newblock {\em Nonconvex semi-linear problems and canonical dual solutions.}
\newblock Advances in Mechanics and Mathematics, Kluwer Academic Publishers, Dordrecht, vol. II, pp. 261–-312 (2003).



\bibitem{gao-opt03}
Gao, D.Y.
\newblock {\em Perfect duality theory and complete solutions to a class of
global optimization problems.}
\newblock Optim., 52(4--5), pp. 467-493 (2003).



\bibitem{gao-cace}
Gao, D. Y.
\newblock {\em Canonical duality theory: theory, method, and
applications in global optimization.}
\newblock Comput. Chem. 33, pp. 1964-1972, (2009).


\bibitem{DGAO2OG}
Gao, D. Y.; Ogden R. W.
\newblock {\em Multiple solutions to non-convex
variational problems with implications for phase transitions and numerical
computation.}
\newblock  Q. J. Mech. Appl. Math, Vol. 61. No. 4, pp. 497-522 (2008).


  \bibitem{gao-ogden08b}Gao, D.Y; Ogden, R.W.
\newblock {\em Closed-form solutions, extremality and nonsmoothness criteria in a large deformation elasticity problem,}
\newblock  Zeitschrift für angewandte Mathematik und Physik (ZAMP),  59, pp.498–517 (2008).

\bibitem{DGAO3RU}
Gao, D. Y.; Ruan N.
\newblock {\em Solutions and optimality criteria for nonconvex quadratic-exponential minimization problem,.}
\newblock Math. Meth. Operations Research, 67 (3), pp. 479-491 (2008).

\bibitem{gao-ruan-pardalos}
Gao, D. Y., Ruan, N; Pardalos, P. M.
\newblock {\em Canonical dual solutions to sum of fourth-order polynomials minimization problems with applications to sensor network localization.}
\newblock Sensors: Theory, Algorithms and Applications, Springer(2010).


\bibitem{gao-sherali09}
Gao, D. Y.; Sherali, H. D.
\newblock {\em Canonical duality: Connection between nonconvex mechanics and global optimization.}
\newblock Advances in Appl. Mathematics and Global Optimization, pp. 249-316, Springer (2009).


\bibitem{gao-gs1}
Gao, D. Y.; Strang, G.
\newblock {\em Geometric nonlinearity: Potential energy, complementary energy, and the gap function.}
\newblock Quart. Appl. Math. 47 (3), pp. 487-504 (1989).


\bibitem{WU}
Gao, D. Y.; Wu, C.
\newblock {\em On the triality theory for a quartic polynomial optimization problem.}
\newblock J. Industrial and Management Optimization, 8(1), pp. 229-242 (2012).

\bibitem{gao-yu}
Gao, D. Y.; Yu, H. F.
\newblock {\em Multi-scale modelling and canonical dual finite element method in phase transitions of solids.}
\newblock Int. J. Solids and Structures, 45, pp. 3660-3673 (2008).



\bibitem{HORN}
Horn, R. A.; Johnson, C. R.
\newblock {\em Matrix Analysis.}
\newblock Cambridge University Press (1985).

\bibitem{marsden-hugh}
Marsden, J. E.; Hughes, T. J. R.
\newblock {\em Mathematical Foundations of Elasticity.}
\newblock Prentice-Hall, (1983).


\bibitem{MAXIMA} Maxima.sourceforge.net.
\newblock {\em Maxima, a Computer Algebra System.}
\newblock Version 5.22.1 (2010). http://maxima.sourceforge.net/


\bibitem{more-wu}
Mor\'{e}, J.; Wu, Z.
\newblock {\em Global continuation for distance geometry problems.}
\newblock SIAM Journal on Optimization, 7, pp. 814-836 (1997).


\bibitem {moreau}
Moreau, J. J.
\newblock {\em La notion de sur-potentiel et les liaisons unilat\'{e}rales en \'{e}lastostatique.}
\newblock C.R. Acad. Sc. Paris, 267 A, pp. 954-957 (1968).


\bibitem {more-pana-stra-88}
Moreau, J. J.; Panagiotopoulos, P. D.; Strang, G.
\newblock {\em Topics in nonsmooth mechanics.}
\newblock Birkhuser Verlag, Basel-Boston, MA. (1988).


\bibitem{PETERS}
Peters, G.; Wilkinson J. H.
\newblock {\em The least squares problem and pseudo-inverses.}
\newblock The Computer Journal, Vol. 13, No 3, pp. 309-316 (1970).

\bibitem{ruan-gao-ima}
Ruan, N.; Gao, D. Y.
\newblock {\em Canonical duality approach for non-linear dynamical systems.}
\newblock IMA J. Appl. Math,  79(2),  313 - 325 (2014).

\bibitem{ruan-gao-pe} Ruan, N.; Gao, D.Y.
\newblock {\em Global optimal solutions to a general
sensor network localization problem,}
\newblock  { Performance Evaluations} 75-76: 1–16 (2014).

\bibitem{r-g-j}
Ruan, N.; Gao, D. Y.; Jiao, Y.
\newblock {\em Canonical dual least square method for solving general nonlinear systems of quadratic equations.}
\newblock Computational Optimization and Applications, Vol 47, pp. 335-347 (2010).


\bibitem{santos-gao}
Santos, H. A. F. A.; Gao D. Y.
\newblock {\em Canonical dual finite element method for solving post-buckling problems of a large deformation elastic beam.}
\newblock Int. J. Nonlinear Mechanics, 47, 240-247 (2012).  doi:10.1016/j.ijnonlinmec.2011.05.012 (2011).


\bibitem{saxe}
Saxe, J.
\newblock {\em Embeddability of weighted graphs in k-space is strongly NP-hard.}
\newblock Proc. 17th Allerton Conference in Communications, Control, and Computing, Monticello, IL, pp. 480-489 (1979).



\bibitem{SEWELL}
Sewell, M. J.
\newblock {\em Maximum and minimum principles.}
\newblock Cambridge University Press, Cambridge, New York, Port
Chester, Melbourne Sydney (1987).


\bibitem{v-z} Strugariu, R.; Voisei, M.D.;  Zalinescu, C.
\newblock {\em Counter-examples in bi-duality, triality
and tri-duality,}
\newblock { Discrete and Continuous  Dynamical Systems}, 31, 1453 - 1468 (2011).

\bibitem{vz2} Voisei, M.D.; Zalinescu, C.: 
\newblock  {\em Some remarks concerning Gao-Strang's
 complementary gap function,}
 \newblock { Applicable Analysis}, {  90}, 1111-1121 (2010).

\bibitem{wang-etal}
Wang, Z. B.; Fang, S. C.; Gao, D. Y.; Xing, W. X.
\newblock {\em Canonical dual approach to solving the maximum cut problem.}
\newblock   { J. Glob. Optim.}, {  54}, 341--352 (2012).


\bibitem{xie}
Xie D.; Schlick T.
\newblock {\em Visualization of chemical databases using the singular
value decomposition and truncated-Newton minimization.}
\newblock Optimization in Computational Chemistry and Molecular Biology, Kluwer Academic Publishers B.V. (2000).


\bibitem{zgy}
Zhang J.; Gao, D. Y.; Yearwood, J.
\newblock {\em A novel canonical dual computational approach for prion AGAAAAGA amyloid fibril molecular modeling.}
\newblock Journal of Theoretical Biology, 284, pp. 149-157 (2011). doi:10.1016/j.jtbi.2011.06.024


\bibitem{ZIA}
Zia, R. K. P.; Redish, E. F.; McKay S. R.
\newblock {\em Making Sense of the Legendre Transform}
\newblock American Journal of Physics, Vol. 77, Issue 7, pp. 614-622
(2009).


\end{thebibliography}
\end{document}